%% file: Cographic-toricARXIVFINAL.tex
\documentclass[11 pt]{amsart}

\usepackage{amssymb,amsmath,amsthm,amscd,mathrsfs,graphicx, color}
\usepackage[matrix,arrow,tips,curve]{xy}
\usepackage[english]{babel}
\usepackage[leqno]{amsmath}
\usepackage{amssymb,amsthm}
\usepackage{amscd}
\usepackage{enumerate}
\usepackage{layout}

\newtheorem{thm}{Theorem}[section]

\newtheorem{lemma}[thm]{Lemma}

\newtheorem{prop}[thm]{Proposition}

\newtheorem{cor}[thm]{Corollary}

\newtheorem{prob}[thm]{Problem}

\newtheorem{teoalpha}{Theorem}

\theoremstyle{remark}
\newtheorem{remark}[thm]{Remark}

\theoremstyle{definition}
\newtheorem{defi}[thm]{Definition}
\newtheorem{Section: Notation}[thm]{}
\newtheorem{exa}[thm]{Example}

\numberwithin{equation}{section}

\newenvironment{sis}{\left\{\begin{aligned}}{\end{aligned}\right.}

\newcommand{\ov}{\overline}
\newcommand{\un}{\underline}

\newcommand{\Supp}{\operatorname{Supp}}

\newcommand{\Spec}{\operatorname{Spec}}

\newcommand{\supp}{\operatorname{Supp}}

\newcommand{\rk}{\operatorname{rk}}

\newcommand{\calA}{{ \mathcal A}}
\newcommand{\calB}{{ \mathcal B}}
\newcommand{\calC}{{ \mathcal C}}

\newcommand{\calF}{{ \mathcal F}}

\newcommand{\calH}{{ \mathcal H}}

\newcommand{\calL}{{ \mathcal L}}

\newcommand{\calN}{{ \mathcal N}}
\newcommand{\calO}{{ \mathcal O}}

\newcommand{\bbG}{{\mathbb G}}

\newcommand{\bbN}{{\mathbb N}}

\newcommand{\bbQ}{{\mathbb Q}}
\newcommand{\bbR}{{\mathbb R}}

\newcommand{\bbZ}{{\mathbb Z}}

\newcommand{\p}{\mathfrak p}
\newcommand{\m}{\mathfrak{m}}
\newcommand{\n}{\mathfrak{n}}

\newcommand{\Gm}{\bbG_m}

\newcommand{\norFan}[1][\vor_{\Gamma}]{\calN(#1)}

\newcommand{\vor}{\operatorname{Vor}}

\renewcommand{\OE}{\overset{\rightarrow}{E}}
\newcommand{\OS}{\overset{\rightarrow}{S}}
\newcommand{\el}[1][e]{\overset{\leftarrow}{#1}}
\newcommand{\er}[1][e]{\overset{\rightarrow}{#1}}

\newcommand{\OC}{\overset{\rightarrow}{C}}
\newcommand{\Cyc}{\overset{\rightarrow}{{\rm Cir}}}
\newcommand{\Cir}{{\rm Cir}}

\newcommand{\OP}{\mathcal{OP}}

\newcommand{\Str}{\rm Str}

\newcommand{\SimCom}{\Delta(\Cyc(\Gamma))}

\newcommand{\vol}{\rm vol}

\begin{document}

\title[Cographic rings]{The geometry and combinatorics of\\  cographic toric face rings}

\author[Casalaina-Martin]{Sebastian Casalaina-Martin}
\address{University of Colorado at Boulder, Department of Mathematics,
Campus Box 395, Boulder, CO-80309-0395, USA.}
\email{casa@math.colorado.edu}

\author[Kass]{Jesse Leo Kass}
\address{University of Michigan, Department of Mathematics, 530 Church St, Ann
Arbor, MI 48103, USA}
\email{jkass@umich.edu}

\author[Viviani]{Filippo Viviani}
\address{Dipartimento di Matematica,
Universit\`a Roma Tre,
Largo S. Leonardo Mur-ialdo 1,
00146 Roma (Italy)}
\email{viviani@mat.uniroma3.it}
\thanks{The first author was partially supported by the NSF grant DMS-1101333.  The second author was partially supported by the NSF grant RTG DMS-0502170.  The third author is a member of
the research center CMUC (University of Coimbra) and is supported by the FCT project \textit{Espa\c cos de Moduli em Geometria Alg\'ebrica} (PTDC/MAT/111332/2009) and by PRIN project
\textit{Geometria delle variet\`a algebriche e dei loro spazi di moduli} (funded by MIUR, cofin 2008).}

\date{\today}



\begin{abstract}
In this paper we define and study a ring associated to a graph that we call the cographic toric face ring, or simply the cographic ring.  The cographic ring is the toric face ring defined by the following equivalent combinatorial structures of a graph: the cographic arrangement of hyperplanes, the Voronoi polytope, and the poset of totally cyclic orientations.  We describe the properties of the cographic ring and, in particular, relate the invariants of the ring to the invariants of the corresponding graph.

Our study of the cographic ring fits into a body of work on describing rings constructed from graphs.  Among the rings that can be constructed from a graph, cographic rings are particularly interesting because they appear in the study of compactified Jacobians of nodal curves.
\end{abstract}

\maketitle

\bibliographystyle{amsalpha}

\section*{Introduction}\label{Sec: Intro}
In this paper we define and study a ring $R(\Gamma)$ associated to a graph $\Gamma$ that we call the cographic toric face ring, or simply the cographic ring.  The cographic ring $R(\Gamma)$ is the toric face ring defined by the following equivalent combinatorial structures of $\Gamma$: the cographic arrangement of hyperplanes $\calC_{\Gamma}^{\perp}$, the Voronoi polytope $\vor_{\Gamma}$, and the poset of totally cyclic orientations $\OP_{\Gamma}$.   We describe the properties of the cographic ring and, in particular, relate the invariants of the ring to the invariants of the corresponding graph.

Our study of the cographic ring fits into a body of work on describing rings constructed from graphs.  Among the rings that can be constructed from a graph, cographic rings are particularly interesting because they appear in the study of compactified Jacobians.

The authors establish the connection between $R(\Gamma)$ and the local geometry of compactified Jacobians in  \cite{CMKV}.  The compactified Jacobian $\bar{J}_{X}^{d}$ of a nodal curve $X$ is the coarse moduli space parameterizing sheaves on $X$ that are rank $1$, semi-stable, and of fixed degree $d$.  These moduli spaces have been constructed by Oda--Seshadri \cite{OS}, Caporaso \cite{caporaso},  Simpson \cite{simpson}, and Pandharipande \cite{Pan}, and the different constructions are reviewed in \cite[Section~2]{CMKV}.    In \cite[Theorem~A]{CMKV}, the authors proved that the completed local ring of $\bar{J}_{X}^{d}$ at a point is isomorphic to a power series ring over  the completion of $R(\Gamma)$, for a graph $\Gamma$ constructed from the dual graph of $X$.

In \cite{CMKV}, the authors also study the local structure of the universal compactified Jacobian, which is a family of varieties over the moduli space of stable curves whose fibers are closely related the compactified Jacobians just discussed.  (See \cite[Section~2]{CMKV} for a discussion of the relation between the compactified Jacobians from the previous paragraph and the fibers of the universal Jacobian).  Caporaso first constructed the universal compactified Jacobian in \cite{caporaso}, and Pandharipande gave an alternative construction in \cite{Pan}.   In \cite[Theorem~A]{CMKV}, the authors  gave a presentation of the completed local ring of the universal compactified Jacobian at a point, and they will explore the relation between that ring and the affine semigroup ring defined in Section~\ref{Ss:aff-cog} in the upcoming paper \cite{CMKV2}.

Cographic toric face rings are examples of toric face rings.  Recall that a toric face ring is constructed from the same combinatorial data that is used to construct a toric variety: a fan.  Let $H_{\bbZ}$ be a a free, finite rank $\bbZ$-module and $\mathcal F$  be a  fan that decomposes
$H_{\mathbb R}=H_{\mathbb Z} \otimes_{\mathbb Z}\mathbb R$  into (strongly convex rational polyhedral) cones.  Consider the free $k$-vector space with basis given by monomials $X^c$ indexed by elements $c \in H_{\bbZ}$.  If we define a multiplication law on this vector space by setting
$$X^c\cdot X^{c'}=
\begin{sis}
& X^{c+c'} & \text{ if  $c, c' \in \sigma\text{ for some } \sigma \in \mathcal F$,} \\
& 0 &\text{ otherwise, }
\end{sis}
$$
and extending by linearity, then the resulting ring $R(\calF)$ is the toric face ring (over $k$) that is associated to $\calF$.

We define the cographic toric face ring $R(\Gamma)$ of a graph $\Gamma$ to be toric face ring associated to the fan that is defined by the cographic arrangement $\calC_{\Gamma}^{\perp}$.  The cographic arrangement is an arrangement of hyperplanes in the real vector space $H_{\bbR}$ associated to the homology group $H_{\bbZ} := H_{1}(\Gamma, \bbZ)$ of the graph.  Every  edge of $\Gamma$ naturally induces a functional on $H_{\bbR}$, and the zero locus of this functional is a hyperplane in $H_{\bbR}$, provided the functional is nonzero.  The cographic arrangement is defined to be the collection of all hyperplanes constructed in this manner.  The intersections of these hyperplanes define a fan $\calF_{\Gamma}^{\perp}$, the cographic fan.  The toric face ring associated to this fan is $R(\Gamma)$.

We study the fan $\calF_{\Gamma}^{\perp}$ in Section~\ref{Sec: comparison}.  The main result of that section is Corollary~\ref{Cor: DelPoset}, which provides two alternative descriptions of $\calF_{\Gamma}^{\perp}$.  First, using a theorem of Amini, we prove that $\calF_{\Gamma}^{\perp}$ is equal to the normal fan of the Voronoi polytope $\vor_{\Gamma}$.  As a consequence, we can conclude that $\calF_{\Gamma}^{\perp}$, considered as a poset, is isomorphic to the poset of faces of $\vor_{\Gamma}$ ordered by reverse inclusion.  Using work of Greene--Zaslavsky, we show that this common poset is also isomorphic to the poset $\OP_{\Gamma}$ of totally cyclic orientations.

The combinatorial definition of $R(\Gamma)$ does not appear in \cite{CMKV}.  Rather,
the rings in that paper appear as invariants under a torus action.  The following theorem, proven in Section~\ref{Sec: presen} (Theorem~\ref{Thm: Pres-inv}), shows that the rings in \cite{CMKV} are (completed) cographic rings.

\begin{teoalpha}\label{teob}  Let $\Gamma$ be a finite  graph with vertices $V(\Gamma)$, oriented edges $\OE(\Gamma)$ and source and target maps $s,t:\OE(\Gamma)\to V(\Gamma)$.
Let
$$
T_\Gamma:=\prod_{v\in V(\Gamma)}\mathbb G_m \ \ \text { and } \ \ A(\Gamma) :=  \frac{k[U_{\el}, U_{\er}\: :\: e\in E(\Gamma)]} {(U_{\el} U_{\er} \colon e \in E(\Gamma))}.
$$
If we make $T_{\Gamma}$ act on $A(\Gamma)$ by
$$
	\lambda \cdot U_{\er}= \lambda_{s(\er)}\:U_{\er}\: \lambda_{t(\er)}^{-1},
$$
then the invariant subring $A(\Gamma)^{T_\Gamma}$ is isomorphic to the cographic ring $R(\Gamma)$.
\end{teoalpha}

The cographic ring $R(\Gamma)$ has reasonable geometric properties.  Specifically, in Theorem~\ref{T:prop-cog} we prove that $R(\Gamma)$ is
\begin{itemize}
	\item of pure dimension $b_1(\Gamma)=\dim_{\mathbb R}H_1(\Gamma,\mathbb R)$;
	\item Gorenstein;
	\item semi-normal;
	\item semi log canonical.
\end{itemize}

We also compute  invariants of $R(\Gamma)$ in terms of the combinatorics of $\Gamma$.
The invariants we compute are the following:
\begin{itemize}
	\item a description of $R(\Gamma)$ in terms of oriented subgraphs (see Section~\ref{Ss:desc-cog});
	\item the number of minimal primes in terms of orientations (Theorem~\ref{T:prop-cog}\eqref{I:minpr-cog});
	\item the embedded dimension of $R(\Gamma)$ in terms of circuits  (Theorem~\ref{T:prop-cog}\eqref{I:emb-cog});
	\item the multiplicity of $R(\Gamma)$ (Theorem~\ref{T:prop-cog}\eqref{I:mult-cog}).
\end{itemize}
Finally, it is natural to ask what information is lost in passing from $\Gamma$ to $R(\Gamma)$.  An answer to this question
is given by Theorem~\ref{Thm: Tor}, which states that $R(\Gamma)$ determines $\Gamma$ up to $3$-edge connectivization.

Combinatorially defined rings, such as the cographic toric face ring, have long been used in the study of
compactified Jacobians, and, more generally, degenerate abelian varieties (e.g., \cite{mumford72}, \cite{OS}, \cite{faltings}, \cite{nam80}, \cite{AN}, \cite{alex}).  In particular, the ring $R(\Gamma)$ we study here is a special case of the rings $R_0(c)$ studied
by Alexeev and Nakamura \cite[Theorem~3.17]{AN}.  There the rings appear naturally as a by-product of  Mumford's technique for degenerating an abelian variety.
Alexeev and Nakamura  proved that $R_0(c)$ satisfies the Gorenstein condition in \cite[Lemma~4.1]{AN}, and the semi-normality was established by Alexeev in  \cite{alexeev02}.  In personal correspondence, Alexeev informed  the  authors that the techniques of those papers can also be used to establish other results in this paper, such as the fact that $R(\Gamma)$ is semi log canonical.

In a different direction, the cographic ring is defined by the cographic fan $\calF_{\Gamma}^{\perp}$, which is the normal fan to the Voronoi polytope $\vor_{\Gamma}$.  There is a body of work studying similar polytopes and the algebro-geometric objects defined by these polytopes.  In \cite{AH}, Altmann--Hille define the polytope of flows associated to an oriented graph (or quiver).  Associated to this polytope is a toric variety which they relate to a moduli space. There are also a number of recent papers that  study the modular/integral flow polytope in $H_{1}(\Gamma, \bbR)$.  This study is motivated by the work of Beck--Zaslavsky on interpreting graph polynomials in terms of lattice points \cite{BZ}.  Some recent papers on this topic are \cite{BZ}, \cite{dall10}, \cite{breuer12}, and \cite{chen10}.  The paper \cite{dall10}, in particular, studies graph polynomials using tools from commutative algebra.   The Voronoi polytope does not equal the modular/integral flow polytope or the polytope of flows of an oriented graph.  It would, however, be interesting to further explore the relation between these polytopes.  (We thank the anonymous referee for pointing out this literature.)

This paper suggests several other questions for further study. First, in Section~\ref{Ss:aff-cog} we exhibit a collection of generators $V_{\gamma}$, indexed by oriented circuits $\gamma$, for $R(\Gamma \setminus T, \phi)$.  What is an explicit set of generators for the ideal of relations between the variables $V_{\gamma}$?  This problem is posed as Problem~\ref{Prob: Circuits}.  Second, in Theorem~\ref{T:prop-cog} we give a formula for the multiplicity of $R(\Gamma)$ in terms of the subdiagram volume of certain semigroups associated to $\Gamma$.  Problem~\ref{Prob: Multiplicity} is: find an expression for this multiplicity in terms of well-known graph theory invariants.  Third, we also prove in Theorem~\ref{T:prop-cog} that $\Spec( R(\Gamma))$ is semi log canonical. In Problem~\ref{Rmk: sdlt}, we ask: which graphs $\Gamma$ have the stronger property that $R(\Gamma)$ is semi divisorial log canonical?

\subsection*{Acknowledgements}
This work began when the authors were visiting the MSRI, in Berkeley, for the special semester in Algebraic Geometry in the spring of 2009; we would like
to thank the organizers of the program as well as the institute for the excellent working conditions and the stimulating atmosphere.
We would like to thank Bernd Sturmfels for his interest in this work, especially for some useful suggestions regarding multiplicity, and for pointing out a mistake in a previous version of this paper.
We thank Farbod Shokrieh for some comments on an early draft of this manuscript and for pointing out the connection between some of our results and the theory of oriented matroids.  We thank the referees for many useful comments  and for pointing out a mistake in a previous version of this paper.

\section{Preliminaries} \label{Sec: prelim}

In this section we review the definitions of the  graph-theoretic objects considered in this paper.   This will provide the reader with enough background to follow the main ideas of the proof of Theorem~\ref{teob} (proven in Section~\ref{Sec: presen}), as well as the proofs of many of the geometric properties of cographic rings (proven in Section~\ref{Sec: Toricfan}).

\subsection{Notation}

Following Serre's notation in \cite[Section~2.1]{SerT}, a \emph{graph} $\Gamma$ will consist of the data
$(\OE  \xymatrix{ \ar @{->}^s @< 2pt> [r] \ar@{->}_t @<-2pt> [r] &  } V, \OE \stackrel{\iota}{\to} \OE),$
where $V$ and $\OE$ are sets, $\iota$ is a fixed-point free involution, and $s$ and $t$ are maps satisfying $s(\er)=t(\iota(\er))$ for all $\er \in \OE$.  The maps $s$ and $t$ are called the \emph{source} and \emph{target} maps respectively.  We call $V=:V(\Gamma)$ the set of \emph{vertices}.  We call $\OE=:\OE(\Gamma)$ the set of \emph{oriented edges}.  We define the set of \emph{(unoriented) edges} to be  $E(\Gamma)=E:=\OE/\iota$.    An \emph{orientation of an edge} $e\in E$ is a representative for $e$ in $\OE$; we use the notation $\er$ and $\el$ for the two possible orientations of $e$.   An \emph{orientation of a graph $\Gamma$} is a section $\phi:E\to \OE$ of the quotient map.  An \emph{oriented graph}  consists of a pair $(\Gamma,\phi)$ where $\Gamma$ is a graph and $\phi$ is an orientation.    Given an oriented graph, we say that $\phi(e)$ is the \emph{positive orientation} of the edge $e$. Given a subset $S\subseteq E$, we define $\OS\subseteq \OE$ to be the set of all orientations of the edges in $S$.

\subsection{Homology of a Graph}\label{SubSec: HomGra}

Given a ring $A$, let $C_0(\Gamma,A)=\OC_0(\Gamma,A)$ be the free $A$-module with basis $V(\Gamma)$ and
$\OC_1(\Gamma,A)$ be the $A$-module generated by
$\OE(\Gamma)$ with the relations $\el=-\er$ for every $e\in E(\Gamma)$.  If we fix an orientation,
then a basis for $\OC_1(\Gamma, A)$ is given by the positively oriented edges; this induces an isomorphism with the usual group of $1$-chains on the simplicial complex associated to $\Gamma$.
These modules may be put into a chain complex.  Define a boundary map $\partial$ by
\begin{gather*}
	\partial: \OC_1(\Gamma,A) \to \OC_0(\Gamma,A)=C_0(\Gamma,A) \\
	\er  \mapsto t(\er)-s(\er).
\end{gather*}
We will denote by $H_\bullet (\Gamma,A)$ the groups obtained from the homology of $\OC_\bullet (\Gamma,A)$.
The  homology groups $H_\bullet(\Gamma,A)$ coincide with the homology groups of the topological space associated to $\Gamma$.

\subsection{The bilinear form}\label{SubSec: scalprod}
The vector space $\OC_1(\Gamma,\mathbb R)$ is endowed with a positive definite bilinear form
$$
(\ , \ ):\OC_1(\Gamma,\mathbb R)\otimes \OC_1(\Gamma,\mathbb R)\to \mathbb R.
$$
uniquely determined by $(\er,\er)=1$, $(\er,\el)=-1$ and $(\er, \er[f])=0$ if $ \er[f]\ne \er, \el$. As above, fixing an orientation induces a basis for $\OC_1(\Gamma,\mathbb R)$, and in terms of such a basis this is the standard inner product.  By restriction, we get a positive definite bilinear form on $H_1(\Gamma,\bbR)\subseteq \OC_1(\Gamma,\mathbb R)$. The pairing $(\cdot, \cdot)$ allows us to form the product $(\er, v)$ of an oriented edge $\er$ with a vector $v\in \OC_1(\Gamma,\mathbb R)$, but not the product $(e, v)$ of $v$ with an unoriented vector.
However, we will write $(e,v)=0$ to mean $(\er, v)=0$ for one (equivalently all) orientations of $e$.

\subsection{Cographic arrangement}\label{SubSec: cographic}
We review the definition of the \emph{cographic arrangement} $\calC_{\Gamma}^{\perp}$ of $\Gamma$ (see \cite[Section~8]{GZ} and \cite[Section~5]{NPS}).\footnote{The name cographic arrangement suggests the fact  that $\calC_{\Gamma}^{\perp}$ depends on the cographic matroid associated to $\Gamma$
The notation $\calC_{\Gamma}^{\perp}$ is used in \cite{NPS}, while in \cite{GZ} the cographic arrangement is denoted by $\calH^{\perp}[\Gamma]$.
There is a dual notion, namely that of the graphic arrangement, which depends only on the graphic matroid associated to $\Gamma$ and is denoted by
$\calC_{\Gamma}$ in \cite{NPS} and by $\calH[\Gamma]$ in \cite[Section~7]{GZ}.
The graphic arrangement of hyperplanes is also studied in \cite[Section~2.4]{OT},
where it is denoted by $\calA(\Gamma)$.}
To begin, let $\mathscr  H$
be the coordinate hyperplane arrangement in $\OC_1(\Gamma,\mathbb R)$.  More precisely:
$$
\mathscr H = \bigcup_{e\in E} \{v\in \OC_1(\Gamma,\mathbb R): (v,e)=0\}.
$$
The restriction of this hyperplane arrangement to $H_1(\Gamma,\mathbb R)$
is called the cographic arrangement $\calC_{\Gamma}^{\perp}$.
 More precisely
$$
\calC_{\Gamma}^{\perp} = \bigcup_{e\in E, H_1(\Gamma,\mathbb R)\nsubseteq \ker (\cdot , e)} \{v\in \OC_1(\Gamma,\mathbb R): (v,e)=0\}.
$$

The cographic arrangement partitions $H_{1}(\Gamma, \bbR)$ into a finite collection of strongly convex rational polyhedral cones.  These cones, together with their faces, form
a (complete) fan that is defined to be the \emph{cographic fan}, and is denoted $\calF_{\Gamma}^{\perp}$.\footnote{We use the notation $\calF_{\Gamma}^{\perp}$ and the name cographic fan  in order to be coherent with the notation $\calC_{\Gamma}^{\perp}$ used in \cite{NPS}
for the cographic arrangement of hyperplanes.}
We give a more detailed enumeration of the cones of this fan later in Section~\ref{Sec: comparison} when we discuss the poset of totally cyclic orientations.

\begin{remark}\label{remcycleedge}
The following observation used in the proof of  Theorem~\ref{teob} is proven in Corollary~\ref{Cor: cyccone}.  We emphasize it here so that the reader may follow the proof of Theorem~\ref{teob} having read just Section~\ref{Sec: prelim}.
\emph{Let  $c=\sum_{e\in E}a_e \er$ and $c'=\sum_{e\in E}a_e'\er$ be cycles in $H_1(\Gamma,\mathbb Z)$.  Then $c$ and $c'$ lie in a common cone of $\calF_{\Gamma}^{\perp}$ if and only if for all $e\in E$,  $a_ea_e'\ge 0$.}  In words, two cycles lie in a common cone if and only if every common edge is oriented in the same direction.
\end{remark}

\subsection{Toric face rings}

We recall the definition of a toric face ring associated to a fan. In \cite[Section~2]{IR} and \cite[Section~2]{BKR}, the authors define more generally the toric face ring associated to a monoidal complex.  The following definition is a special case.

\begin{defi}\label{dfntfr}
Let  $H_{\mathbb Z}$  a free $\mathbb Z$-module of finite rank and
let  $\mathcal F$ be a fan of  (strongly convex rational polyhedral) cones in  $H_{\mathbb R}=H_{\mathbb Z} \otimes_{\mathbb Z}\mathbb R$ with support $\supp \calF$.
 The  \emph{toric face ring}  $R_k(\mathcal F)$   is the  $k$-algebra whose underlying $k$-vector space
has basis $\{X^c\: :\: c\in H_{\mathbb Z}\cap \supp \calF\}$ and whose multiplication is defined by
\begin{equation} \label{Eqn: CographicTimesLaw}
X^c\cdot X^{c'}=
\begin{sis}
& X^{c+c'} & \text{ if } c, c' \in \sigma\text{ for some } \sigma \in \mathcal F \\
& 0 &\text{ otherwise. }
\end{sis}
\end{equation}
We will write $R(\calF)$ if we do not need to specify the base field $k$.
\end{defi}

\begin{remark}\label{remred}
It follows from the  definition that $R(\mathcal F)$ is a reduced ring finitely generated over $k$.  See also Section~\ref{S:prop-cog}, especially \eqref{E:pres2-cog}, for more on generators and relations.
\end{remark}

A cographic toric face ring is a toric face ring associated to a cographic fan.

\begin{defi}\label{dfnctfr} Let $\Gamma$ be a finite graph.
The \emph{cographic toric face ring}  $R_k(\Gamma)$ is the toric face $k$-ring $R(\calF_{\Gamma}^{\perp})$
associated to the cographic fan $\calF_{\Gamma}^{\perp}$.  We will write $R(\Gamma)$ if we do not need to specify the base field $k$.
\end{defi}

\subsection{The Voronoi polytope}
Following \cite{BDN}, we define the \emph{Voronoi polytope} of $\Gamma$ by
\begin{displaymath}
	\vor_{\Gamma} := \{v \in H_{1}(\Gamma, \bbR) : (v,v) \le (v-\lambda, v-\lambda) \text{ for all } \lambda \in H_1(\Gamma,\bbZ)  \}.
\end{displaymath}
The reader familiar with the Voronoi decomposition of $\bbR^n$ will recognize this polytope as the unique cell containing the origin in the Voronoi decomposition associated with the lattice $H_1(\Gamma,\bbZ)$ endowed with the scalar product defined in \S\ref{SubSec: scalprod} (see \cite{Erd} or \cite[Section~2.5]{alex} for more details).

To the Voronoi polytope, we can associate its \emph{normal fan} $\mathcal N(\vor_{\Gamma})$ which is defined as follows.  Given a face $\delta$ of $\vor_{\Gamma}$, we define the (strongly convex rational polyhedral) cone $C_{\delta}$ by
$$
C_\delta=
\{\alpha\in H_1(\Gamma,\mathbb R): (\alpha,r) \ge (\alpha,r') \ \ \text{ for all } \ \ r\in \delta, r'\in \vor_{\Gamma}\}.
$$
The normal fan $\mathcal N(\vor_{\Gamma})$ of $\vor_{\Gamma}$ is the fan whose cones are the cones $C_{\delta}$.

\begin{remark} \label{remfanpoly}
In Proposition~\ref{Prop:  FanCompare} we will prove that the cographic fan $\calF_{\Gamma}^{\perp}$  is equal to the normal fan  of the Voronoi polytope $\mathcal N(\vor_{\Gamma})$.
\end{remark}

\section{Totally cyclic orientations} \label{Sec: orient}
Here we define and study totally cyclic orientations of a graph.  We also define an oriented circuit on a graph and describe the relation between these circuits and totally cyclic orientations.

\subsection{Subgraphs}
 \label{SubSec: Del-Contr}
In this subsection, we introduce some special subgraphs that will play an important role  throughout the paper.

Given a graph $\Gamma$ and a collection $S \subset E(\Gamma)$ of edges, we define $\Gamma \setminus S$
to be the graph, called a  \emph{spanning subgraph} (e.g., see \cite[\S4]{OS}),  obtained from $\Gamma$ by removing the edges in $S$ and leaving the vertices unmodified.
In other words $\Gamma \setminus S$ consists of the data
$$(\overset{\rightarrow}{E(\Gamma)\setminus S}  \xymatrix{ \ar @{->}^s @< 2pt> [r] \ar@{->}_t @<-2pt> [r] &  } V, \overset{\rightarrow}{E(\Gamma)\setminus S} \stackrel{\iota}{\to} \overset{\rightarrow}{E(\Gamma)\setminus S}).$$
Of particular significance is the special case where $S = \{e\}$ consists of a single edge.  If $\Gamma \setminus \{ e \}$ has more connected components than $\Gamma$, then we say that $e$ is a \emph{separating edge}.
The \emph{set of all separating edges} is written $E(\Gamma)_{\rm sep}$.

Given a chain $c\in \vec C_1(\Gamma,\mathbb R)$ we would like to refer to the underlying graph having only those edges in the support of  $c$.  More precisely, given $c\in \vec C_1(\Gamma,\mathbb R)$, let $\Supp(c)$ denote the set of all edges $e$ with the property that $(e,c) \ne 0$.  We define $\Gamma_{c}$ to be the subgraph of $\Gamma$ with $V(\Gamma_{c}) := V(\Gamma)$ and $E(\Gamma_{c}) := \Supp(c)$.  There is a distinguished orientation $\phi_{c}$ of $\Gamma_{c}$ given
by setting $\phi_{c}(e)$ equal to  $\er$ if $(\er,c) >0$ and to $\el$ otherwise.
Using this subgraph, we can write $c$ as
\begin{equation}\label{Eqn: CanForm} 	
c=\sum_{e\in \Supp(c)} m_{c}(e) \phi_{c}(e),
\end{equation}
with all $m_{c}(e) >0$.  Indeed, we have $m_{c}(e) = (\phi_c(e), c)$.

\subsection{Totally cyclic orientations and oriented circuits}\label{SubSec: Circ}

Totally cyclic orientations will play a dominant role in what follows. We are going to review their definition and their basic properties.

\begin{defi}  \label{Def: TotCyclic}
	If $\Gamma$ is connected, then we say that an orientation $\phi$ of $\Gamma$ is {\it totally cyclic} if there does not exist a non-empty proper subset $W\subset V(\Gamma)$ such that every edge $e$ between a vertex in $W$ and a vertex in the complement $V(\Gamma)\smallsetminus W$ is oriented from $W$ to $V \smallsetminus W$ (i.e.~the source of $\phi(e)$ lies in $W$ and the target of $\phi(e)$ lies in $V(\Gamma) \smallsetminus W$).
	If $\Gamma$ is disconnected, then we say that an orientation of $\Gamma$ is totally cyclic if the orientation induced on each connected component of $\Gamma$ is totally cyclic.
\end{defi}

Observe that, if $\Gamma$ is a graph with no edges, then the empty orientation of $\Gamma$ is a totally cyclic orientation.
Totally cyclic orientations are closely related to oriented circuits.  Recall that a graph $\Delta$ is called  \emph{cyclic} if it is connected, free from separating edges,
and satisfies $b_1(\Delta)=1$. A cyclic graph together with a totally cyclic orientation is called an \emph{oriented circuit}.
  A cyclic graph admits exactly two totally cyclic orientations.

Let $\Cyc(\Gamma)$ denote the set of all oriented circuits on $\Gamma$; that is, $\gamma=(\Delta,\phi_{\Delta})$ is an element of $\Cyc(\Gamma)$ if
$\Delta$ is a cyclic subgraph of  $\Gamma$ and $\phi_{\Delta}$ is a totally cyclic orientation of $\Delta$. We call $E(\Delta)$ the \emph{support} of
$\gamma=(\Delta,\phi_{\Delta})\in \Cyc(\Gamma)$. There is a natural map
$$
\Cyc(\Gamma)\to H_1(\Gamma,A)
$$
given by
$$
\gamma=(\Delta,\phi_\Delta) \mapsto [\gamma]= \sum_{e\in E(\Delta)} \phi_\Delta(e).
$$

With respect to the orientation $\phi$ of $\Gamma$, we can consider the subset $\operatorname{Cir}_{\phi}(\Gamma)\subset \Cyc(\Gamma)$ that
consists of oriented circuits on $\Gamma$  of the form $(\Delta, \phi|_{\Delta})$ (i.e.~oriented circuits  whose orientation is compatible with $\phi$).

\begin{remark}\label{Rmk: orien-mat}
The set $\Cyc(\Gamma)$ of oriented circuits on $\Gamma$ are the (signed) cocircuits
of the cographic oriented matroid $M^*(\Gamma)$ of $\Gamma$ or, equivalently, the (signed) circuits of the oriented
graphic matroid $M(\Gamma)$ of $\Gamma$ (see \cite[Section~1.1]{ormat}).  Many of the combinatorial results that follow can be naturally stated using this language.
We will limit ourselves to pointing out the connection with the theory, when relevant.
\end{remark}

The next lemma clarifies the relationship between totally cyclic orientations and compatibly oriented circuits.
Recall that an oriented path from $w\in V(\Gamma)$ to $v\in V(\Gamma)$ is a collection of oriented edges $\{\er_1,\ldots,\er_r\}\subset \OE(\Gamma)$
such that $s(\er_1)=w$, $t(\er_i)=s(\er_{i+1})$ for any $i=1,\ldots r-1$, and $t(\er_r)=v$. If $\phi$ is an orientation of $\Gamma$, a path
compatibly oriented with respect to $\phi$ is an oriented path as before of the form $\{\phi(e_1),\ldots,\phi(e_r)\}$.

\begin{lemma}\label{Lemma: TotCyclic}
Let $\Gamma$ be a  graph.
\begin{enumerate}[(1)]
	\item \label{Item: TotCycExist}
			The graph $\Gamma$ admits a totally cyclic orientation if and only if $E(\Gamma)_{\rm sep}=\emptyset$.
	 \item \label{Item: IsTotCycO}

		Fix an orientation $\phi$ on $\Gamma$. The following conditions are equivalent:
		\begin{enumerate}[(a)]
			\item \label{Item: TotCyc}
				The orientation   is totally cyclic.
			\item \label{Item: vTow}
				 For any distinct $v, w\in V(\Gamma)$ belonging to the same connected component of $\Gamma$, there exists a path compatibly oriented with respect to $\phi$ from $w$ to $v$.
			\item \label{Item: Gener}
				 The cycles $[\gamma]$ associated to the  $\gamma \in \operatorname{Cir}_{\phi}(\Gamma)$ generate $H_1(\Gamma,\bbZ)$, and $E(\Gamma)_{\rm sep}=\emptyset$.
			\item \label{Item: EdgeToCyc}
				Every edge $e \in E$ is contained in the support of a compatibly oriented circuit $\gamma \in \operatorname{Cir}_{\phi}(\Gamma)$.
		\end{enumerate}
\end{enumerate}
\end{lemma}
\begin{proof}
For part~\eqref{Item: TotCycExist}, see e.g.~\cite[Lemma~2.4.3(1)]{CV1} and the references therein.
Part~\eqref{Item: IsTotCycO} is a reformulation of \cite[Lemma~2.4.3(2)]{CV1}. The only difference is that part~\eqref{Item: IsTotCycO} is proved in loc.~cit.~under the additional hypothesis that
$E(\Gamma)_{\rm sep}=\emptyset$.
Note however that each of the conditions (a), (b) and (d) imply that $E(\Gamma)_{\rm sep}=\emptyset$, hence we deduce part~\eqref{Item: IsTotCycO} as stated above.
\end{proof}

The following well-known lemma can be thought of as a modification of (c) above.  We no longer require that the oriented circuits on $\Gamma$ be oriented compatibly.   The statement is essentially that any cycle $c$ in $H_1(\Gamma,\mathbb Z)$ is a positive linear combination of cycles associated to circuits supported on $c$.

\begin{lemma}\label{Lemma: DecCyc} Let $\Gamma$ be a graph
	and let $c \in \OC_1(\Gamma,\bbZ)$.
	Then $c \in H_1(\Gamma,\bbZ)$ if and only if $c$ can be expressed as
	\begin{equation} \label{Eqn: CycleSum}
		c=\sum_{\gamma\in \operatorname{Cir}_{\phi_c}(\Gamma_c)} n_{c}(\gamma) [\gamma],
	\end{equation}
	for some natural numbers $n_{c}(\gamma)\in \bbN$.
\end{lemma}
\begin{proof}
A direct proof follows from the definitions and is left to the reader. Alternatively, one can use the fact that a covector of an oriented matroid can be written as a composition of cocircuits conformal to it (see \cite[Proposition~3.7.2]{ormat}) together with Remark~\ref{Rmk: orien-mat}.
\end{proof}

The oriented circuits can be used to define a  simplicial complex which will be used in Section~\ref{Ss:desc-cog}.

\begin{defi}\label{Def: Concor}
	Two oriented circuits $\gamma=(\Delta,\phi)$ and $\gamma'=(\Delta',\phi')$ are said to be
\emph{concordant}, written $\gamma \asymp \gamma'$, if for any $e\in E(\Delta)\cap E(\Delta')$ we have $\phi(e)=\phi'(e)$.
We write $\gamma\not\asymp\gamma'$ if $\gamma$ and $\gamma'$ are not concordant.
\end{defi}

\begin{defi}\label{Def: SimplComp}
	The \emph{simplicial complex of concordant circuits}, $\SimCom$, is defined to be the (abstract) simplicial complex whose elements are collections $\sigma \subseteq \Cyc(\Gamma)$ of oriented
	circuits on $\Gamma$ with the property that any two circuits are concordant (i.e.~if $\gamma_1,\gamma_2\in \sigma$, then $\gamma_1 \asymp \gamma_2$).
\end{defi}

\subsection{The poset $\OP_{\Gamma}$ of totally cyclic orientations} \label{SubSec: TotOrient}

Totally cyclic orientations naturally form a poset.
This poset was defined in \cite[Section~5.2]{CV1}, but we recall the definition for the sake of completeness.

\begin{defi}\cite[Definition~5.2.1]{CV1}  \label{Def: TotCycOrient}
	The poset $\OP_{\Gamma}$ of \emph{totally cyclic orientations} of $\Gamma$ is the set of pairs
	$(T, \phi)$ where $T\subset E(\Gamma)$ and
	$\phi:E(\Gamma\setminus T) \to \OE(\Gamma\setminus T)$
	is a totally cyclic orientation of
	$\Gamma\smallsetminus T$,\footnote{The choice of orientation on the complement of $T$, rather than on $T$ itself, has to do with the importance of the notion of spanning subgraphs of $\Gamma$, all of which are of this form.  In graph theory it is customary to denote spanning
subgraphs in this way, so we follow that convention.}  endowed with the following partial order
	\begin{displaymath}
		(T',\phi')  \le (T, \phi) \Leftrightarrow   \Gamma \setminus T'\subseteq  \Gamma \setminus T
		\text{ and } \phi'=\phi|_{E(\Gamma\smallsetminus  T')}.
	\end{displaymath}
	We call $T$ the {\it support} of the pair $(T,\phi)$.
\end{defi}
Using Lemma~\ref{Lemma: TotCyclic}\eqref{Item: EdgeToCyc}, we get that
\begin{equation} \label{Eqn: IneqCyc}
	(T',\phi')\le (T,\phi) \Leftrightarrow \Cir_{\phi'}(\Gamma \setminus T')\subseteq \Cir_{\phi}(\Gamma \setminus T).
\end{equation}
The set $\Cir_{\phi}(\Gamma \setminus T)$ is a collection of concordant cycles.
Another connection between orientations and totally cyclic orientations is given by the following
definition.

\begin{defi}
	Let $\sigma \in \SimCom$ be a collection of concordant circuits.  To $\sigma$ we associate
	the pair $(T_{\sigma}, \phi_{\sigma}) \in \OP_{\Gamma}$ which is defined as follows.  Set $T_{\sigma}$
	equal to the set of all edges that are \emph{not}  contained in a circuit $\gamma \in \sigma$.
	The orientation $\phi_{\sigma}$ of $\Gamma \setminus T_{\sigma}$ is defined by setting
$$\phi_{\sigma}(e):=
\begin{cases}
\er & \text{ if } (\er,[\gamma])>0 \text{ for all } \gamma \in \sigma, \\
\el & \text{ if } (\el,[\gamma])>0 \text{ for all } \gamma \in \sigma. \\
\end{cases}$$
\end{defi}

Observe that the orientation $\phi_{\sigma}$ on $\Gamma\setminus T_{\sigma}$ is a totally cyclic orientation by  Lemma~\ref{Lemma: TotCyclic}\eqref{Item: EdgeToCyc}
and that $\sigma\subseteq \Cir_{\phi_{\sigma}}(\Gamma\setminus T_{\sigma})$.
The following lemma will be useful in the sequel.

\begin{lemma}\label{L:max-orie}
The maximal elements of the poset $\OP_{\Gamma}$ are given by $(E(\Gamma)_{\rm sep}, \phi)$ as $\phi$ varies among the totally cyclic	orientations of $\Gamma \setminus E(\Gamma)_{\rm sep}$.
\end{lemma}
\begin{proof}
The proof is   left to the reader.
\end{proof}

\begin{remark}\label{Rmk: covectors}
The poset $\OP_{\Gamma}$ of totally cyclic orientations is isomorphic to the poset of covectors of the cographic oriented matroid $M^*(\Gamma)$ of $\Gamma$ (see \cite[Section~3.7]{ormat}). Equivalently, the poset obtained from $\OP_{\Gamma}$
by adding an element $\un{1}$ and declaring that $\un{1}\geq (T,\phi)$ for any $(T,\phi)\in \OP_{\Gamma}$
is isomorphic to the big face lattice $\calF_{\rm big}(M^*(\Gamma))$ of the cographic oriented matroid $M^*(\Gamma)$
(see \cite[Section~4.1]{ormat}).

\end{remark}

\section{Comparing posets: The cographic arrangement, the Voronoi polytope and totally cyclic orientations} \label{Sec: comparison}

In this section we prove that the poset $\OP_{\Gamma}$ of totally cyclic orientations of $\Gamma$ is isomorphic to the poset of cones (ordered by inclusion) of  the cographic fan $\calF_{\Gamma}^{\perp}$, which we also show is  the normal fan of the Voronoi polytope $\vor_{\Gamma}$ of $\Gamma$.

\subsection{Cographic arrangement}\label{SubSec: arrang}
Let us start by describing the cographic arrangement $\calC_{\Gamma}^{\perp}$ associated to  $\Gamma$  in the language of totally cyclic orientations.

For every edge $e\in E(\Gamma)$, we can consider the linear subspace of $H_1(\Gamma,\bbR)$
\begin{displaymath} 	
\{(\cdot, e)=0\}:=\{ v \in H_{1}(\Gamma, \bbR) : (v,e) =0 \}.
\end{displaymath}
This subspace is a proper subspace (i.e.~a hyperplane) precisely when $e$ is not a separating edge, and the collection of all such hyperplanes is defined to be the cographic arrangement.  Similarly, for any oriented edge $\er\in \OE(\Gamma)$, we set
$$\{(\cdot, \er)\geq 0\}:=\{v\in H_1(\Gamma,\bbR)\: : \: (v, \er)\geq 0\}.$$

As mentioned above, the elements of the cographic arrangement partition $H_{1}(\Gamma, \bbR)$ into a finite collection of rational polyhedral cones.  These cones, together with their faces, form
the cographic fan $\calF_{\Gamma}^{\perp}$.  We can enumerate these cones, and make their relation to totally cyclic orientations more explicit by introducing some notation.

Given a collection $T$ of edges and
an orientation $\phi$ of $\Gamma \setminus T$ (not necessarily totally cyclic), we define (possibly empty) cones $\sigma(T, \phi)$ and $\sigma^{\operatorname{o}}(T, \phi)$ by
\begin{align}
	\sigma(T,\phi):=& \bigcap_{e\not\in T}\{(\cdot, \phi(e)) \ge 0\} \ \cap \
		\bigcap_{e\in T}\{(\cdot, e)=0\},   \label{Eqn: FacDel} \\
		\sigma^{\operatorname{o}}(T,\phi):=& \bigcap_{e\not\in T}\{(\cdot, \phi(e)) > 0\} \ \cap \
		\bigcap_{e\in T}\{(\cdot, e)=0\}.  \label{Eqn: intcones}
\end{align}

The cone $\sigma^{\operatorname{o}}(T, \phi)$ is a subcone of $\sigma(T, \phi)$, and it  is the relative interior of $\sigma(T, \phi)$ provided
$\sigma^{\operatorname{o}}(T, \phi)$  is non-empty.
The cone $\sigma(T, \phi)$ is an element of the cographic fan, and every cone in the fan can be written in this form.
While every element of $\calF_{\Gamma}^{\perp}$ can be written as $\sigma(T, \phi)$,
the pair $(T, \phi)$ is \emph{not} uniquely determined by the cone.  The pair $(T, \phi)$ is, however, uniquely determined if we further require that $(T, \phi) \in \OP_{\Gamma}$.  This fact is proven in the following proposition, which is essentially a restatement of some results of Greene--Zaslavsky (\cite[Section~8]{GZ}).

\begin{prop}\label{Prop: PosetArr}
\noindent
\begin{enumerate}[(i)]
\item \label{Item: Arr1} Every cone $\sigma \in  \calF_{\Gamma}^{\perp}$ can be written as $\sigma = \sigma(T, \phi)$ for
a unique element $(T, \phi) \in \OP_{\Gamma}$.

\item \label{Item: Arr2} For any  $(T,\phi)\in \OP_{\Gamma}$, the linear span of $\sigma(T,\phi)$ is equal to
$$\langle \sigma(T,\phi) \rangle=\bigcap_{e\in T}\{(\cdot, e)=0\}=H_1(\Gamma\setminus T,\bbR)$$
and has dimension $b_1(\Gamma\setminus T)$.
\item \label{Item: Arr3}
For any  $(T,\phi)\in \OP_{\Gamma}$,
the extremal rays of $\sigma(T,\phi)$ are the rays generated by the elements $[\gamma]$ for
$\gamma\in \Cir_{\phi}(\Gamma\setminus T)$.
\end{enumerate}
\end{prop}

\begin{proof}
Part~\eqref{Item: Arr1} follows from \cite[Lemma~8.2]{GZ}. Note that in loc.~cit.~the authors assume that $E(\Gamma)_{\rm sep}=\emptyset$. However, it is easily checked
that the inclusion map $\Gamma\setminus E(\Gamma)_{\rm sep}\subseteq \Gamma$ induces natural isomorphisms
$\calF_{\Gamma\setminus E(\Gamma)_{\rm sep}}^{\perp}\cong \calF_{\Gamma}^{\perp}$ and $\OP_{\Gamma\setminus E(\Gamma)_{\rm sep}}\cong \OP_{\Gamma}$.
Therefore, the general case follows from the special case treated in loc.~cit.

Let us now prove part~\eqref{Item: Arr2}.
The linear subspace $\displaystyle \bigcap_{e\in T}\{(\cdot, e)=0\}\subseteq H_1(\Gamma,\bbR)$ is generated by all the cycles of $\Gamma$ that do not contain edges $e\in T$
in their support and is therefore equal to  $H_1(\Gamma\setminus T,\bbR)$, which has dimension equal to $b_1(\Gamma\setminus T)$.
Now, to complete the proof, let us establish that $\langle \sigma(T,\phi) \rangle=\bigcap_{e\in T}\{(\cdot, e)=0\}$.  First, if  $\sigma(T,\phi)^\circ = \emptyset$, i.e. if
$\sigma(T,\phi)= \{0\}$, then $b_1(\Gamma\setminus T)=0$ by Lemma \ref{Lemma: TotCyclic}\eqref{Item: EdgeToCyc}.
   But then $\bigcap_{e\in T}\{(\cdot, e)=0\}=H_1(\Gamma\setminus T,\bbR)=0$, and we are done.
On the other hand if $\sigma(T,\phi)^\circ \ne \emptyset$,
then  $\sigma^{\operatorname{o}}(T,\phi)$
is the relative interior of $\sigma(T,\phi)$, and hence the linear span of $\sigma(T,\phi)$ is equal to $\displaystyle \bigcap_{e\in T}\{(\cdot, e)=0\}$.

Finally, let us prove part~\eqref{Item: Arr3}. From \cite[Lemma~8.5]{GZ}, it follows that
 the extremal rays of $\sigma(T,\phi)$ are among the rays generated by the elements $[\gamma]$ for $\gamma\in \Cir_{\phi}(\Gamma\setminus T)$. We conclude by showing
that for any $\gamma\in \Cir_{\phi}(\Gamma\setminus T)$, the ray generated by $[\gamma]$ is extremal for $\sigma(T,\phi)$. By contradiction, suppose that we can write
\begin{equation}\label{Eqn: Ext-ray}
[\gamma]=\sum_{\stackrel{\gamma'\in \Cir_{\phi}(\Gamma\setminus T)}{\gamma'\neq \gamma}}m_{\gamma'}[\gamma'],
\end{equation}
for some $m_{\gamma'}\in \bbR_{\geq 0}$. Consider a cycle
$\gamma_0\in \Cir_{\phi}(\Gamma\setminus T)\setminus \{\gamma\}$ such that $m_{\gamma_0}>0$ (which clearly exists
since $[\gamma]\neq 0$). Since $\gamma$ and $\gamma_0$ are concordant and distinct, there should exist
an edge $e\in E(\gamma_0)\setminus E(\gamma)$. Now returning to the expression \eqref{Eqn: Ext-ray},  on the left hand
side neither the oriented edge $\er$ nor $\el$ can appear. On the other hand, on the right hand side  the oriented
edge $\phi(e)$ appears with positive multiplicity, because it appears with multiplicity $m_{\gamma_0}>0$
in $m_{\gamma_0}[\gamma_0]$ and all the oriented circuits appearing in the summation are concordant.
This is a contradiction.
\end{proof}

\begin{cor}\label{Cor: PosetArr}
The association
	\begin{equation*}
		\begin{aligned}
			(T,\phi) & \mapsto \sigma(T,\phi)
		\end{aligned}
	\end{equation*}
 defines an isomorphism between the poset of $\OP_{\Gamma}$ and the poset of cones of $\calF_{\Gamma}^{\perp}$ ordered by inclusion.\footnote{Note that the poset of cones of $\calF_{\Gamma}^{\perp}$ is anti-isomorphic to the face poset $\calL(\calC_{\Gamma}^{\perp})$ of the arrangement
$\calC_{\Gamma}^{\perp}$ (see \cite[Definition~2.18]{OT}).}
In particular the number of connected components of the complement of $\calC_{\Gamma}^{\perp}$ in $H_1(\Gamma,\bbR)$ is equal to
the number of totally cyclic orientations on $\Gamma\setminus E(\Gamma)_{\rm sep}$.
\end{cor}
\begin{proof} According to Proposition~\ref{Prop: PosetArr}\eqref{Item: Arr1},
the map in the statement is bijective. We have to show that
$$\sigma(T,\phi)\subseteq \sigma(T',\phi')\Longleftrightarrow (T,\phi)\leq (T',\phi').$$
The implication $\Leftarrow$ is clear by the definition \eqref{Eqn: FacDel} of $\sigma(T,\phi)$.

Conversely, assume that $\sigma(T,\phi)\subseteq \sigma(T',\phi')$.  There is nothing to show if $\sigma(T, \phi) = \{0\}$ is the origin.  Otherwise, by Proposition~\ref{Prop: PosetArr}\eqref{Item: Arr2}, the relative interior $\sigma^{\operatorname{o}}(T,\phi)$
of $\sigma(T,\phi)$ is non-empty, so pick $c \in \sigma^{\operatorname{o}}(T,\phi)$.  By Formula~\eqref{Eqn: intcones}, for every $e\not\in T$ we have that $(c, \phi(e))>0$.
Since $c\in \sigma(T',\phi')$, by Definition~\eqref{Eqn: FacDel}, we must have $e\not\in T'$ and $\phi'(e)=\phi(e)$. This shows that $T\supseteq T'$ and that
$\phi'_{\Gamma\setminus T}=\phi$, or in other words that $(T,\phi)\leq (T',\phi')$.

The last assertion follows from the first one using  the fact that the connected components of the complement
of $\calC_{\Gamma}^{\perp}$ in $H_1(\Gamma,\bbR)$ are the maximal cones in $\calF_{\Gamma}^{\perp}$  and Lemma~\ref{L:max-orie}.
\end{proof}

\begin{remark}
The last assertion of Corollary \ref{Cor: PosetArr} is due to Green--Zaslavsky (see \cite[Lemma~8.1]{GZ}).
Moreover, Greene--Zaslavsky give a formula for the number of totally cyclic orientations
of a graph free from separating edges \cite[Theorem~8.1]{GZ}.
\end{remark}

The following well-known result will play a crucial role in the proof of Theorem~\ref{Thm: Pres-inv}.

\begin{cor}\label{Cor: cyccone}
Let  $$c=\sum_{e\in E}a_e \er \  \ \text{ and } \ \ c'=\sum_{e\in E}a_e'\er$$ be cycles in $H_1(\Gamma,\mathbb Z)$.  Then there is a cone of $\calF_{\Gamma}^{\perp}$ containing  $c$ and $c'$ if and only if for all $e\in E$,  $a_ea_e'\ge 0$.
\end{cor}
\begin{proof}
From Proposition~\ref{Prop: PosetArr}\eqref{Item: Arr1}, it follows that $c$ and $c'$ belong to the same cone of $\calF_{\Gamma}^{\perp}$ if and only if there exists
$(T,\phi)\in \OP_{\Gamma}$ such that $c, c'\in \sigma(T,\phi)$. We conclude by looking at the explicit description \eqref{Eqn: FacDel}.
\end{proof}

\begin{remark}\label{Rmk: mat-hyper}
 Corollary~\ref{Cor: PosetArr}
together with Remark~\ref{Rmk: covectors} imply that the cographic oriented matroid $M^*(\Gamma)$ is represented by
the cographic hyperplane arrangement $\calC_{\Gamma}^{\perp}$, in the sense of \cite[Section~1.2(c)]{ormat}.
Using this, Corollary~\ref{Cor: cyccone} is a restatement of the fact that two elements of $H_1(\Gamma,\bbZ)$ belong
to the same cone of $\calF_{\Gamma}^{\perp}$ if and only if their associated covectors are conformal
(see \cite[Section~3.7]{ormat}).
\end{remark}

\subsection{Voronoi polytope}\label{SubSec: Vor}

The following description of the faces of $\vor_{\Gamma}$ is a restatement, in our notation, of a result
of Omid Amini (see \cite{Ami}), which gives a positive answer to a conjecture of Caporaso--Viviani (see \cite[Conjecture~5.2.8(i)]{CV1}).

\begin{prop}[\bf{Amini}]
\label{Prop: PosetFaces}
\noindent
\begin{enumerate}[(i)]
\item \label{Item: Fac1} Every face of the Voronoi polytope $\vor_{\Gamma}$ is of the form
\begin{equation}\label{Eqn: FacVor}
		F(T,\phi):=\left\{v\in \vor_{\Gamma}\: : \: (v, [\gamma])=\frac{1}{2}([\gamma],[\gamma]) \text{ for any }
        \gamma\in \Cir_{\phi}(\Gamma\setminus T)\right\}
\end{equation}
for some uniquely determined element $(T,\phi)\in \OP_{\Gamma}$.
\item \label{Item: Fac2} For any  $(T,\phi)\in \OP_{\Gamma}$, the dimension of the affine span
of $F(T,\phi)$ is equal to $b_1(\Gamma(T))=b_1(\Gamma)-b_1(\Gamma\setminus T)$.
\item \label{Item: Fac3} For any  $(T,\phi)\in \OP_{\Gamma}$, the codimension one faces of $\vor_{\Gamma}$ containing $F(T,\phi)$
are exactly those of the form $F(S,\psi)$ where $(S,\psi)\leq (T,\phi)$ and $b_1(\Gamma\setminus S)=1$.
\end{enumerate}
\end{prop}
\begin{proof}
Part~\eqref{Item: Fac1} follows by combining \cite[Theorem~1]{Ami} and \cite[Lemma~7]{Ami}.
Part~\eqref{Item: Fac2} follows from the remark after \cite[Lemma~10]{Ami}.
Part~\eqref{Item: Fac3} follows from \cite[Lemma~7]{Ami}.
\end{proof}

\begin{cor}[\bf{Amini}]
\label{Cor: PosetFaces}
The association
	\begin{equation*}
		\begin{aligned}
			(T,\phi) & \mapsto F(T,\phi)
		\end{aligned}
	\end{equation*}
defines an isomorphism of posets between the poset $\OP_{\Gamma}$ and the poset of faces of $\vor_{\Gamma}$ ordered by reverse inclusion.
In particular the number of vertices of $\vor_{\Gamma}$ is equal to
the number of totally cyclic orientations on $\Gamma\setminus E(\Gamma)_{\rm sep}$.
\end{cor}
\begin{proof}
The first statement is a reformulation of \cite[Theorem~1]{Ami}.
The last assertion follows form the first one together with Lemma~\ref{L:max-orie}.
\end{proof}

We now show that the cographic fan $\calF_{\Gamma}^{\perp}$ is the \emph{normal fan} $\norFan$ of the Voronoi polytope
$\vor_{\Gamma}$.
The cones of the normal fan, ordered by inclusion, form a poset that is clearly isomorphic   to the poset of faces of $\vor_{\Gamma}$, ordered by reverse inclusion.

\begin{prop} \label{Prop: FanCompare}
    The cographic fan $\calF_{\Gamma}^{\perp}$ is equal to  $\norFan$, the normal fan of the Voronoi polytope
    $\vor_{\Gamma}$.
\end{prop}
\begin{proof}
By Propositions~\ref{Prop: PosetArr} and \ref{Prop: PosetFaces}, it is enough to show that, for any $(T,\phi)\in \OP_{\Gamma}$, the normal cone in $\norFan$
to  the face $F(T,\phi)\subset \vor_{\Gamma}$ is equal to $\sigma(T,\phi)$.
Fix a face $F(T,\phi)$ of $\vor_{\Gamma}$ for some $(T,\phi)\in \OP_{\Gamma}$.
If $(T,\phi)$ is equal to the minimal element $\un{0}=(E(\Gamma)_{\rm sep},\emptyset)$ of the poset $\OP_{\Gamma}$
then $F(\un{0})=\vor_{\Gamma}$ and its normal cone is equal to the origin in $H_1(\Gamma,\bbR)$ which is equal to $\sigma(\un{0})$.

 Suppose now that $b_1(\Gamma\setminus T)\geq 1$. Denote by  $\{(S_i,\psi_i)\}$ all the elements of $\OP_{\Gamma}$ such
that $(S_i,\psi_i)\leq (T,\phi)$ and $b_1(\Gamma\setminus S_i)=1$. Let $\gamma_i$ be the unique oriented circuit of $\Gamma$ such that
$\Cir_{\psi_i}(\Gamma\setminus S_i)=\{\gamma_i\}$.
According to Proposition~\ref{Prop: PosetFaces}\eqref{Item: Fac3}, the codimension
one faces of $\vor_{\Gamma}$ containing $F(T,\phi)$ are exactly those of the form $F(S_i,\psi_i)$. Therefore the normal cone of $F(T,\phi)$ is the cone whose extremal rays
are the normal cones to the faces $F(S_i,\psi_i)$ which, using \eqref{Eqn: FacVor}, are equal to $\sigma(S_i,\psi_i)=\bbR_{\geq 0}\cdot [\gamma_i]$.  By Proposition~\ref{Prop: PosetArr}
\eqref{Item: Arr3}, the cone whose extremal rays are given by $\bbR_{\geq 0}\cdot [\gamma_i]$ is equal to $\sigma(T,\phi)$, which completes the proof.
\end{proof}

Combining Corollary~\ref{Cor: PosetArr}, Corollary~\ref{Cor: PosetFaces} and Proposition~\ref{Prop: FanCompare},
we get the following incarnations of the poset $\OP_{\Gamma}$ of totally cyclic orientations.

\begin{cor}\label{Cor: DelPoset}
	The following posets are isomorphic:
	\begin{enumerate}
		\item the poset $\OP_{\Gamma}$ of totally cyclic orientations;
		\item the poset of faces of the Voronoi polytope $\vor_{\Gamma}$, ordered by reverse inclusion;
		\item the poset of cones in the normal fan $\norFan$, ordered by inclusion;
		\item the poset of cones in the cographic fan $\calF_{\Gamma}^{\perp}$, ordered by inclusion.
	\end{enumerate}
\end{cor}

\begin{remark}\label{Rmk: mat-poly}
Corollary~\ref{Cor: DelPoset}
together with Remark~\ref{Rmk: covectors} imply that the cographic oriented matroid $M^*(\Gamma)$ is represented by
the Voronoi polytope $\vor_{\Gamma}$ (which is a zonotope, see e.g.~\cite{Erd}), in the sense of \cite[Section~2.2]{ormat}.
\end{remark}

\section{Geometry of toric face rings} \label{Sec: Toricfan}

Let  $H_{\mathbb Z}$  be a free $\mathbb Z$-module of finite rank $b$ and
let  $\mathcal F$ be a fan of  (strongly convex rational polyhedral) cones in  $H_{\mathbb R}=H_{\mathbb Z} \otimes_{\mathbb Z}\mathbb R$.
The aim of this section is to study the toric face ring $R(\calF)=R_k(\calF)$ associated to $\calF$ as in Definition~\ref{dfntfr}.
We will pay special attention to fans $\calF$ that are {\em complete}, i.e.~such that every $x\in H_{\bbR}$ is contained in some cone $\sigma\in\calF$,
or {\em polytopal}, i.e.~the normal fans of rational poytopes in $H_{\bbR}^*$. Note that a polytopal fan is complete, but the converse is false if $b\geq 3$ (see \cite[page~84]{Oda} for
an example).
In the subsequent sections, we will apply the results of this section to the cographic fan
$\calF_{\Gamma}^{\perp}$ of a graph $\Gamma$, which is polytopal by Proposition~\ref{Prop: FanCompare}.

Note that the fan $\calF$ is naturally a poset: given $\sigma,\sigma'\in \calF$, we say that $\sigma\geq \sigma'$ if $\sigma \supseteq \sigma'$.
The poset $(\calF,\geq)$ has some nice properties, which we now describe.
  Recall the following standard concepts from poset theory.
 A (finite) poset $(P,\leq)$ is called a
\emph{meet-semilattice} if every two elements $x, y\in P$ have a meet (i.e.~an
element, denoted by $x\wedge y$ that is uniquely
characterized by conditions  $x\wedge y\leq x, y$ and, if $z\in P$ is such that
$z\leq x, y$, then $z\leq x\wedge y$).  In a meet-semilattice
every finite subset of elements $\{x_1, \cdots, x_n\}\subset P$ admits a meet, denoted by
$x_1\wedge\cdots\wedge x_n$.
A meet-semilattice is called
\emph{bounded} (from below) if it has a minimum element $\un{0}$.
A bounded meet-semilattice is called \emph{graded} if, for every element
$x\in P$, all  maximal chains from $\un{0}$ to $x$ have the same length.
If this is the case, we define a function, called the \emph{rank function},
$\rho:P\to \bbN$ by setting $\rho(x)$  equal to the length of any maximal chain
from $\un{0}$ to $x$.
A graded meet-semilattice is said to be \emph{pure} if all the maximal
elements have the same rank, and this maximal rank is called the \emph{rank} of the
poset and it is denoted by $\rk P$. A graded meet-semilattice is said to be \emph{generated in maximal rank}
if every element of $P$ can be obtained as the meet of a subset consisting
of maximal elements.

Having made these preliminary remarks, we now collect
some of the properties of the poset $(\calF,\geq)$ which we will need later.

\begin{lemma}\label{L:Fposet}
	The poset $(\calF,\geq)$ has the following properties.
	\noindent
	\begin{enumerate}[(i)]
		\item\label{I:Fposet1}
			$(\calF, \geq)$ is a meet-semilattice, where the meet of two cones is equal to their intersection.
		\item\label{I:Fposet2}
			 $(\calF,\geq)$ is bounded with minimum element $\un{0}$ given by the zero cone $\{0\}$.
		\item\label{I:Fposet3}
			 $(\calF,\geq)$ is a graded semi-lattice with
			rank function given by  $\rho(\sigma):=\dim \sigma$.
		\item\label{I:Fposet4}
			 If $\calF$ is complete, then $(\calF,\geq)$  is pure of rank $\rk \calF=b$.
		\item\label{I:Fposet5}
			 If $\calF$ is complete, then $(\calF,\geq)$ is generated in maximal rank.
	\end{enumerate}
\end{lemma}
\begin{proof}
The proof is left to the reader.
\end{proof}

We will denote by $\calF_{\rm max}$ the subset of $\calF$ consisting of the {\em maximal cones} of $\calF$.

\subsection{Descriptions of $R(\calF)$ as an inverse limit and as a quotient.}\label{Ss:2desc}

In this subsection, we give two descriptions of the toric face ring $R(\calF)$.

The first description of $R(\mathcal F)$ is as  an inverse limit of affine semigroup rings.
For any cone $\sigma\in \calF$, consider  the semigroup
\begin{equation}\label{E:smgr}
C(\sigma):=\sigma\cap H_{\bbZ}\subset H_{\bbZ}.
\end{equation}
which, according to Gordan's Lemma (e.g.~\cite[Proposition~6.1.2]{BH}),  is a positive normal affine semigroup, i.e.~a finitely generated semigroup isomorphic to a subsemigroup of $\bbZ^d$ for some
$d\in \bbN$ such that $0$ is the unique invertible element and such that if $m\cdot z\in C(\sigma)$ for some $m\in \bbN$ and $z\in \bbZ^d$ then $z\in C(\sigma)$.

\begin{defi}\label{D:Aff-smgr}
We define $R_k(\sigma):=k[C(\sigma)]$ to be the affine semigroup ring associated to $C(\sigma)$
(in the sense of \cite[Section~6.1]{BH}), i.e.~the $k$-algebra whose underlying vector space has basis $\{X^c\: :\: c\in C(\sigma)\}$
and whose multiplication is defined by $X^c\cdot X^{c'}:=X^{c+c'}$.
We will write $R(\sigma)$ if we do not need to specify the base field $k$.   If $\mathcal F_\sigma$ is the fan induced by $\sigma$ (consisting of the cones in $\mathcal F$ that are faces of $\sigma$), then clearly $R(\sigma)=R(\mathcal F_\sigma)$.
\end{defi}

The following properties are well-known

\begin{lemma}\label{L:pro-affsmgr}
$R(\sigma)$ is a normal, Cohen--Macaulay domain of dimension equal to $\dim \sigma$.
\end{lemma}
\begin{proof}
By definition, we have  $R(\sigma) \subset k[H_{\bbZ}]=k[x_1^{\pm 1},\dots,x_{b}^{\pm 1}]$, hence $R(\sigma)$ is a domain.
$R(\sigma)$ is normal by \cite[Theorem~6.1.4]{BH} and Cohen--Macaulay
by a theorem of Hochster (see \cite[Theorem~6.3.5(a)]{BH}).
Finally, it follows easily from \cite[Proposition~6.1.1]{BH} that the (Krull) dimension of $R(\sigma)$ is equal to $\dim  \sigma$.
\end{proof}

Given two elements $\sigma, \sigma' \in \calF$ such that $\sigma\geq \sigma'$,
or equivalently such that $\sigma \supseteq \sigma'$, there exists a natural projection
map between the corresponding affine semigroups rings of Definition~\ref{D:Aff-smgr}:
$$
\begin{aligned}
r_{\sigma/\sigma'}:R(\sigma)& \twoheadrightarrow R(\sigma')\\
X^c & \mapsto
\begin{cases}
X^c & \text{ if } c\in \sigma' \subseteq \sigma,\\
0 & \text{ if } c\in \sigma \setminus \sigma'.
\end{cases}\\
\end{aligned}
$$
With respect to the above maps, the set  $\{R(\sigma) \: : \: \sigma \in \calF\}$
forms an inverse system of rings. From \cite[Proposition~2.2]{BKR}, we deduce the following description of $R(\calF)$:

\begin{prop}\label{P:pres1}  Let $\mathcal F$ be a fan.
We have an isomorphism
$$R(\calF)=\varprojlim_{\sigma \in \calF } R(\sigma).$$
\end{prop}

\noindent We denote by $r_{\sigma}:R(\calF)\to R(\sigma)$ the natural projection maps.

The second description of $R(\calF)$ is as a quotient of a polynomial ring.
For any cone $\sigma\in \calF$, the semigroup $C(\sigma)=\sigma\cap H_{\bbZ}$ has a unique minimal generating set, called the {\em Hilbert basis} of $C(\sigma)$ and denoted
by $\calH_{\sigma}$ (see \cite[Proposition~7.15]{MS}). Therefore, we have a surjection
\begin{equation}\label{E:sur-smgrp}
\begin{aligned}
\pi_{\sigma}: k[V_{\alpha}\: :\: \alpha \in \calH_{\sigma}] & \twoheadrightarrow R(\sigma) \\
V_{\alpha} & \mapsto X^{\alpha}.
\end{aligned}
\end{equation}
In the terminology of \cite[Chapter~4]{Stu}, the kernel of $\pi_{\sigma}$, which we denote by $I_{\sigma}$,
is the \emph{toric ideal} associated to the subset $\calH_{\sigma}$. In the terminology of \cite[Chapter~II.7]{MS}, $I_{\sigma}$ is the {\em lattice ideal}�
associated with the kernel of the group homomorphism
$$
\begin{aligned}
p_{\sigma}:\bbZ^{\calH_{\sigma}} & \to H_{\bbZ}�\\
\un{u}=\{{u}_{\alpha}\}_{\alpha\in \calH_{\sigma}} & \mapsto \sum_{\alpha\in \calH_{\sigma}} {u}_{\alpha} \alpha.
\end{aligned}
$$
From \cite[Lemma~4.1]{Stu} (see also \cite[Theorem~7.3]{MS}), we get that $I_{\sigma}$ is a binomial ideal with the following explicit presentation
\begin{equation}\label{E:binom-id}
I_{\sigma}=\langle V^{\un{u}}-V^{\un{v}}\: : \: \un{u}, \un{v}\in \bbN^{\calH_{\sigma}}\subset \bbZ^{\calH_{\sigma}} \: \: \text{�with }�\: p_{\sigma}(\un{u})=p_{\sigma}(\un{v})   \rangle,
\end{equation}
where, for any $\un{u}=({u}_{\alpha})_{\alpha\in \calH_{\sigma}}\in \bbN^{\calH_{\sigma}}$, we set
$V^{\un{u}}:=\prod_{\alpha\in \calH_{\sigma}}�V_{\alpha}^{{u}_{\alpha}}\in k[V_{\alpha}\: :\: \alpha \in \calH_{\sigma}]$.

If we set $\displaystyle \calH_{\calF}:=\bigcup_{\sigma\in \calF} \calH_{\sigma}$ then,
from Definition~\ref{dfntfr}, it follows that we have a surjection
\begin{equation}\label{E:surR}
\begin{aligned}
\pi_{\calF}: k[V_{\alpha}\: :\: \alpha\in \calH_{\calF}] & \twoheadrightarrow R(\calF)\\
V_{\alpha} &\mapsto X^{\alpha}.
\end{aligned}
\end{equation}
We denote by $I_{\calF}$ the kernel of $\pi_{\calF}$. In order to describe the ideal $I_{\calF}$,
we introduce the abstract simplicial complex $\Delta_{\calF}$ on the vertex set
$\calH_{\calF}$  whose faces are the collections of elements of $\calH_{\calF}$ that belong to the same cone of $\calF$.
The minimal non-faces of $\Delta_{\calF}$ are formed by pairs $\{\alpha, \alpha'\}$ of elements of $\calH_{\calF}$
such that $\alpha$ and $\alpha'$ do not belong to the same cone of $\calF$; hence $\Delta_{\calF}$ is a flag complex
(see \cite[Chapter~III, Section~4]{Sta}). Consider the  the Stanley--Reisner ring (or face ring)
$$k[\Delta_{\calF}]:=\frac{k[V_{\alpha}\: :\: \alpha\in \calH_{\calF}]} {(V_{\alpha}V_{\alpha'}\: :\: \{\alpha,\alpha'\}\not\in \Delta_{\calF})}$$
associated to the flag complex $\Delta_{\calF}$
(see \cite[Chapter~II]{Sta} for an introduction to Stanley--Reisner rings).
Observe that if $\{\alpha,\alpha'\}\not\in \Delta_{\calF}$ then $X^{[\alpha]}\cdot X^{[\alpha']}=0$ by Definition~\ref{dfntfr}. This implies that the surjection $\pi_{\calF}$ factors as
\begin{equation*}
\pi_{\calF}: k[V_{\alpha}\: :\: \alpha\in \calH_{\calF}]  \twoheadrightarrow
\frac{k[V_{\alpha}\: :\: \alpha\in \calH_{\calF}]}{(V_{\alpha} V_{\alpha'}\: :\: \{\alpha,\alpha'\}\not\in \Delta_{\calF})} =k[\Delta_{\calF}]\twoheadrightarrow R(\calF),
\end{equation*}
or in other words that $(V_{\alpha}V_{\alpha'}\: :\: \{\alpha,\alpha'\}\not\in \Delta_{\calF})\subseteq I_{\calF}$.

Moreover, observe also that the surjection $\pi_{\calF}$ of \eqref{E:surR} is compatible with the surjections
$\pi_{\sigma}$ of \eqref{E:sur-smgrp} for every $\sigma\in \calF$, in the sense that
we have a commutative diagram
\begin{equation}\label{E:diarings}
\xymatrix{
k[V_{\alpha}\: :\: \alpha\in \calH_{\calF}]\ar@{->>}[d]_{\theta} \ar@{->>}[r]^(.6){\pi_{\calF}} & R(\calF)
\ar@{->>}[d]^{r_{\sigma}}\\
k[V_{\alpha}\: :\: \alpha\in \calH_{\sigma}] \ar@{->>}[r]^(.6){\pi_{\sigma}}\ar@/_/[u]_s&
R(\sigma)
}
\end{equation}
where $\theta$ is the surjective ring homomorphism given by sending $V_{\alpha}\mapsto V_{\alpha}$ if $\alpha\in \calH_{\sigma}\subseteq \calH_{\calF}$ and to $V_{\alpha}\mapsto 0$ if $\alpha\in \calH_{\calF}\setminus \calH_{\sigma}$.
Both the vertical surjections have natural sections: the left map has a section $s$ obtained by sending
$V_{\alpha}\mapsto V_{\alpha}$ for any $\alpha\in \calH_{\sigma}\subset \calH_{\calF}$ and the left map
has a section obtained by sending $X^{c}$ into $X^{c}$ for any $c\in C(\sigma)=\sigma\cap H_{\bbZ}\subset H_{\bbZ}$.
Therefore we can regard $I_{\sigma}$ as an ideal
of $k[V_{\alpha}\: :\: \alpha\in \calH_{\calF}]$ by extensions of scalars and, by the above commutative diagram,
we have that $I_{\sigma}\subseteq I_{\calF}$.

From \cite[Propositions~2.3 and 2.6]{BKR}, we get the following description of the ideal $I_{\calF}$:

\begin{prop}\label{P:pres2}  Let $\mathcal F$ be a fan.
The kernel $I_{\calF}$ of the map $\pi_{\calF}$ of \eqref{E:surR} is given by
$$
I_{\calF}=(V_{\alpha}V_{\alpha'}\: :\: \{\alpha, \alpha'\}\not\in \Delta_{\calF})+\sum_{\sigma\in \calF}
I_{\sigma}=(V_{\alpha}V_{\alpha'}\: :\: \{\alpha, \alpha'\}\not\in \Delta_{\calF})+\sum_{\sigma\in \calF_{\rm max}}
I_{\sigma},
$$
where, as usual, $\calF_{\rm max}$ denotes the subset of $\calF$ consisting of the maximal cones.
\end{prop}

\subsection{Prime ideals of $R(\calF)$.}\label{Ss:ideals}

We now want to describe the prime ideals of the ring $R(\calF)$.
Observe that, from the Definition~\ref{dfntfr}, it follows that $R(\calF)$ has a natural $\bbZ^{b}\cong H_{\bbZ}$-grading.

Recall the following notions for a $\bbZ^n$-graded ring $R$ (see e.g.~\cite{Uli}).
A \emph{graded ideal} is an ideal $I$ of $R$ with the property that for any
$x \in I$ all homogenous components of $x$ belong to $I$ as well; this is equivalent to
$I$ being generated by homogenous elements. For any ideal $I$ of $R$ the
\emph{graded core} $I^*$ of $I$ is defined as the ideal generated by all homogenous elements of $I$.
It is the largest graded ideal contained in $I$. If $\p$ is a prime ideal of $R$
then $\p^*$ is a prime ideal (see \cite[Lemma 1.1(ii)]{Uli}).

For any $\sigma \in \calF $, the kernel of the natural projection map
$r_{\sigma} :R(\calF)\twoheadrightarrow R(\sigma)$, which is explicitly
equal to
\begin{equation}\label{E:pr-ideals}
{\mathfrak p}_{\sigma}:=(\{X^c\: :\: c\not\in \sigma\}),
\end{equation}
is graded, since it is generated by homogeneous elements, and is prime by Lemma~\ref{L:pro-affsmgr}. From  \cite[Lemma 2.1]{IR}, we deduce the following description of the graded ideals of $R(\calF)$:

\begin{prop}\label{P:pr-ideals}
\noindent
The assignment $\sigma \mapsto {\mathfrak p}_{\sigma}$ gives an isomorphism between the poset $(\calF,\geq)$
and the poset of graded prime ideals of $R(\calF)$ ordered by reverse inclusion.
In particular, $\m=\p_{\{0\}}$ is the unique graded maximal ideal of $R(\calF)$, which is also a maximal ideal in the usual sense.
\end{prop}

From Proposition~\ref{P:pr-ideals}, we can deduce a description of the minimal primes of $R(\calF)$.

\begin{cor}\label{C:min-pr}
The minimal primes of $R(\calF)$ are the primes  $\p_{\sigma}$, as $\sigma$ varies among all the maximal cones of $\calF$.
In particular, if $\calF$ is complete then $R(\calF)$ is of pure dimension $b$.
\end{cor}
\begin{proof}
Observe that if $\p$ is a minimal ideal of $R(\calF)$, then $\p^*=\p$ by the minimality of $\p$; hence $\p$ is graded.
Conversely, if $\p$ is a graded ideal of $R(\calF)$ which is minimal among the graded ideals of $R(\calF)$ then $\p$ is also a minimal ideal of $R(\calF)$: indeed if ${\mathfrak q}\subseteq \p$
then  ${\mathfrak q}^* = \p$ by the minimality properties of $\p$; hence ${\mathfrak q}=\p$.

It is now clear that the first assertion follows from Proposition~\ref{P:pr-ideals}.
The last assertion follows from the first one together with Lemma~\ref{L:Fposet}\eqref{I:Fposet4} and Lemma~\ref{L:pro-affsmgr}.
\end{proof}

\begin{defi}\label{D:strata}
The poset of {\em strata} of $R(\calF)$, denoted by $\Str(R(\calF))$, is the set of all the ideals of $R(\calF)$ that are sums of minimal primes, with the order relation given by reverse inclusion.
\end{defi}

Geometrically, the poset $\Str(R(\calF))$ is the collection of all scheme-theoretic intersections of irreducible components of $\Spec R(\calF)$,
ordered by inclusion.

\begin{cor}\label{C:strata}
If $\calF$ is complete then the assignment $\sigma\mapsto \p_{\sigma}$ gives an isomorphism between $(\calF,\geq)$ and $\Str(R(\calF))$.
\end{cor}
\begin{proof}
The statement will follow from Proposition~\ref{P:pr-ideals} if we show that the ideals that are sums of minimal primes of $R(\calF)$ are exactly those of the form $\p_{\sigma}$, for some $\sigma\in \calF$.
Indeed, given minimal primes $\p_{\sigma_i}$ for $i=1,\ldots,n$ (see Corollary~\ref{C:min-pr}),
we have that $\displaystyle \bigcap_{i=1}^n \sigma_i =\sigma$ for some $\sigma\in \calF$ and, from \eqref{E:pr-ideals}, it follows that
\begin{equation}\label{eqn*}
\sum_{i=1}^n \p_{\sigma_i}=\left(X^c \: : \: c \not\in \bigcap_{i=1}^n \sigma_i \right)=\p_{\sigma}.
\end{equation}
Conversely, every cone $\sigma\in \calF $ is the intersection of the maximal dimensional cones $\sigma_i$
containing it by Lemma~\ref{L:Fposet}\eqref{I:Fposet5}. Therefore \eqref{eqn*} shows that $\p_{\sigma}\in \Str(R(\calF))$.
\end{proof}

\subsection{Gorenstein singularities}\label{Ss:Gor}

The aim of this subsection is to prove the following theorem, that $R(\calF)$ is Gorenstein  provided that $\calF$ is a polytopal fan.

\begin{thm}\label{T:Gore}
If $\calF$ is a polytopal fan, then  $R(\calF)$ is a Gorenstein ring and its canonical module $\omega_{R(\calF)}$ is isomorphic to
$R(\calF)$ as a graded module.
\end{thm}

\begin{proof}
The theorem follows from two results of Ichim--R\"omer \cite{IR}.
 The first is \cite[Theorem~1.1]{IR} stating that a toric  face ring  $R(\mathcal F)$ is Cohen--Macaulay  provided that the fan $\mathcal F$ is \emph{shellable} (see \cite[p.252]{IR} for the definition).
 The second is  \cite[Theorem~1.4]{IR} stating that $R(\mathcal F)$ is Gorenstein and its canonical module $\omega_{R(\mathcal F)}$ is isomorphic to $R(\mathcal F)$ as a graded module provided that $R(\mathcal F)$ is Cohen--Macaulay and $\mathcal F$ is \emph{Eulerian} (see \cite[Definition~6.4]{IR} for the definition).

 Now it is enough to recall that a polytopal fan is Eulerian (see e.g.~\cite[p.302]{Sta3}) and shellable by the
Bruggesser--Manni Theorem (see \cite[Theorem~5.2.14]{BH}).
\end{proof}

\subsection{The normalization}\label{Ss:norm}

In this subsection, we prove that the toric face ring of any fan is seminormal and we describe its normalization.

Recall that, given a reduced ring $R$ with total quotient ring $Q(R)$,
the \emph{normalization} of $R$, denoted by $\ov{R}$, is the integral closure of $R$ inside $Q(R)$.
$R$ is said to be normal if $R=\ov{R}$
(see for example \cite[Definition~1.5.1]{HS}).
Moreover, we need the following

\begin{defi}\label{D:semnorm}
Let $R$ be a Mori ring, i.e.~a reduced ring such that $\ov{R}$ is finite over $R$.
The \emph{seminormalization} of $R$, denoted by ${}^+R$, is the biggest subring of $\ov{R}$ such that
the induced pull-back map $\Spec({}^+R)\to \Spec R$ is bijective with trivial residue field extension.
We say that $R$ is seminormal if ${}^+R=R$.
\end{defi}

For the basic properties of seminormal rings, we refer to \cite{GT} and \cite{Swa}.
Observe that $R(\calF)$ is a Mori ring since it is reduced and finitely generated over a field $k$
(see Remark~\ref{remred}).

\begin{thm}\label{T:semnorm}
Let $\calF$ be any fan.
\noindent
\begin{enumerate}[(i)]
\item \label{I:semnorm1} The normalization of $R(\calF)$ is equal to
$$\ov{R(\calF)}=\prod_{\sigma\in \calF_{\rm max}} R(\sigma),$$
where $\calF_{\rm max}$ is the subset of $\calF$ consisting of all the maximal cones of $\calF$.

\item \label{I:semnorm2} $R(\calF)$ is a seminormal ring.
\end{enumerate}
\end{thm}
\begin{proof}
Let us first prove Part~\eqref{I:semnorm1}. By \cite[Corollary~2.1.13]{HS} and Corollary~\ref{C:min-pr},
we get that the normalization of $R(\calF)$ is equal to
$$\ov{R(\calF)}= \prod_{\sigma\in \calF_{\rm max}} \ov{R(\sigma)}.$$
We conclude by Lemma~\ref{L:pro-affsmgr}, which says that each domain $R(\sigma)$ is normal.

Let us now prove Part \eqref{I:semnorm2}. According to Proposition \ref{P:pres1} and Lemma \ref{L:pro-affsmgr},
the ring $R(\calF)$ is an inverse limit of normal domains. Then the seminormality of $R(\calF)$ follows from \cite[Corollary~3.3]{Swa}.
\end{proof}

\subsection{Semi log canonical singularities}\label{Ss:slc}

In this subsection, we prove that   $\Spec R(\calF)$ has semi log canonical  singularities provided that $\calF$ is a polytopal fan.

We first recall the definitions of  log canonical  and semi log canonical pairs
(see \cite{KM} for log canonical pairs and \cite{AFKM} or \cite{Fuj1} for semi log canonical pairs).
For the relevance of slc singularities in the theory of compactifications of moduli spaces, see \cite{Kol}.

\begin{defi}\label{D:slc}
Let $X$ be an $S_2$ variety (i.e.~such that the local ring $\calO_{X,x}$ of $X$ at any (schematic) point $x\in X$ has depth at least ${\rm min}\{2, \dim \calO_{X,x}\}$) of pure dimension $n$ over a field $k$ and $\Delta$ be an effective $\bbQ$-Weil divisor on $X$ such that $K_X+\Delta$ is $\bbQ$-Cartier.
\noindent
\begin{enumerate}[(i)]
\item \label{I:slc1}We say that the pair $(X,\Delta)$ is log canonical (or \emph{lc} for short)  if
\begin{itemize}
\item $X$ is smooth in codimension one (or equivalently $X$ is normal);
\item There exists a log resolution $f:Y\to X$ of $(X,\Delta)$ such that
$$K_Y=f^*(K_X+\Delta)+\sum_i a_i E_i,$$
where $E_i$ are divisors on $Y$ and $a_i\geq -1$ for every $i$.
 \end{itemize}
 We say that $X$ is lc  if the pair $(X,0)$ is lc, where $0$ is the zero divisor.
\item \label{I:slc2} We say that the pair $(X,\Delta)$ is semi log canonical (or \emph{slc} for short)  if
\begin{itemize}
\item $X$ is nodal in codimension one (or equivalently $X$ is seminormal and Gorenstein in codimension one);
\item If $\mu:X^{\mu}\to X$ is the normalization of $X$ and $\Theta$ is the $\bbQ$-Weil divisor on $X$ given
by
\begin{equation}\label{Eqn: theta}
K_{X^{\mu}}+\Theta=\mu^*(K_X+\Delta),
\end{equation}
then the pair $(X^{\mu}, \Theta)$ is lc.
\end{itemize}
 We say that $X$ is slc  if the pair $(X,0)$ is slc, where $0$ is the zero divisor.
 \end{enumerate}
\end{defi}

With the above definitions, we can prove the following

\begin{thm}\label{T:slc}
If $\calF$ is a polytopal fan then the variety $\Spec R(\calF)$ is slc.
\end{thm}

\begin{proof}
Observe that $\Spec R(\mathcal F)$ is Gorenstein by Theorem \ref{T:Gore} and seminormal
by  Theorem \ref{T:semnorm}\eqref{I:semnorm2}; hence in particular it is $S_2$ and nodal in codimension one (see \cite[Section~8]{GT}). Moreover $\Spec R(\mathcal F)$ is of pure
dimension $\rk \mathcal F$ by  Corollary \ref{C:min-pr}.
Consider now the normalization morphism  (see Theorem \ref{T:semnorm}\eqref{I:semnorm1})
$$\mu:\Spec \ov{R(\mathcal F)}=\coprod_{\sigma\in \calF_{\rm max}} \Spec R(\sigma) \to \Spec R(\mathcal F),$$

If we apply the formula \eqref{Eqn: theta} to the above morphism $\mu$ and we use the fact that $\Delta=0$ (by hypothesis) and $K_{X}=0$ by
Theorem \ref{T:Gore}, then we get that the divisor $\Theta$ restricted to each connected component $ \Spec R(\sigma) $ of the normalization
$\Spec \ov{R(\mathcal F)}$  is equal to $-K_{ \Spec R(\sigma) }$.
Therefore, from Definition \ref{D:slc}\eqref{I:slc2}, we get that
 $\Spec R(\mathcal F)$ is slc if and only if the pair $( \Spec R(\sigma), -K_{ \Spec R(\sigma)})$
is lc for every $\sigma\in \calF_{\rm max}$ . Therefore, we conclude using the fact that for any toric variety $Z$ the pair $(Z,-K_Z)$ is lc (see \cite[Proposition~2.10]{FS} or \cite[Corollary~11.4.25]{CLS}).
\end{proof}

\subsection{Embedded dimension} \label{Ss:embdim}

In this subsection, we compute the embedded dimension of $R(\calF)$ at its unique graded maximal ideal $\m$.
In doing this, we also compute the embedded dimension of the affine semigroup ring $R(\sigma)$ of Definition \ref{D:Aff-smgr} at the maximal ideal $(X^c\: :\: c\in C(\sigma)\setminus \{0\})$ which, by a slight abuse of notation, we also denote by $\m$.

Recall that given a maximal ideal $\m$ of a ring $R$ with residue field $k:=R/\m$, the embedded dimension of $R$ at $\m$ is the dimension of the $k$-vector space
$\m/\m^2$. Geometrically, the embedded dimension of $R$ at $\m$ is the dimension of the Zariski tangent space of $\Spec(R)$ at the point $\m\in \Spec(R)$.

\begin{thm}\label{T:embdim}  Let $\mathcal F$ be a fan.
\noindent
\begin{enumerate}[(i)]
\item \label{I:Emb1} The embedded dimension of $R(\sigma)$ at $\m$ is equal to the cardinality of the Hilbert basis $\calH_{\sigma}$ (Section~\ref{Ss:2desc}).
\item \label{I:Emb2} The embedded dimension of $R(\calF)$ at $\m$ is equal to the cardinality of $\calH_{\calF}$ ($=\bigcup_{\sigma\in \calF} \calH_{\sigma}$).
\end{enumerate}
\end{thm}
\begin{proof}
Consider  the presentation \eqref{E:sur-smgrp} of the ring $R(\sigma)$.
Since the elements of the Hilbert basis $\calH_{\sigma}$ cannot be written in a nontrivial way as $\bbN$-linear combinations of elements in the semigroup $C(\sigma)$ (see the proof of \cite[Prop. 7.15]{MS}), we get that
the ideal $I_{\sigma}=\ker \pi_{\sigma}$ satisfies
\begin{equation}\label{E:incl1}
	 I_{\sigma}\subset \n^2,
\end{equation}
where $\n:=(V_{\alpha}\: :\: \alpha\in \calH_{\sigma})\subset k[V_{\alpha}\: :\: \alpha\in \calH_{\sigma}]$.
Part \eqref{I:Emb1} now follows from \eqref{E:sur-smgrp} and \eqref{E:incl1}.

In order to prove part \eqref{I:Emb2}, consider the presentation \eqref{E:surR} of the ring $R(\calF)$. It is enough to prove that the ideal $I_{\calF}=\ker \pi_{\calF}$ satisfies
\begin{equation}\label{E:incl2}
	 I_{\calF}\subset \mathfrak{o}^2,
\end{equation}
where $\mathfrak{o}:=(V_{\alpha}\: :\: \alpha\in \calH_{\calF})\subset k[V_{\alpha}\: :\:
\alpha\in \calH_{\calF}]$. Consider the generators of $I_{\calF}$ given in Proposition \ref{P:pres2}.
Clearly the generators of the form $V_{\alpha}V_{\alpha'}$ (for $\{\alpha,\alpha'\}\not\in \Delta_{\calF}$) belong to $\mathfrak{o}^2$.
In order to deal with the other generators of $I_{\calF}$, consider the diagram \eqref{E:diarings}. As in the discussion that precedes Proposition
\ref{P:pres2}, we view $I_{\sigma}$ as included in $I_{\calF}$ via the section $s$.
By applying the section $s$ to the inclusion \eqref{E:incl1} and using the obvious inclusion $s(\n^2)\subseteq \mathfrak{o}^2$, we get the desired inclusion \eqref{E:incl2}.
\end{proof}

\subsection{Multiplicity} \label{Ss:mult}

In this subsection, we study the multiplicity $e_{\m}(R(\calF))$ of $R(\calF)$ at its unique graded maximal
ideal $\m$.

Recall (see for example \cite[Chapter~IIB, Theorem~3]{Ser}) that the  Hilbert--Samuel function
$$n\mapsto \dim_k R(\calF)/\m^n,  $$
is given, for large values of $n\in \bbN$, by a polynomial (called the Hilbert--Samuel polynomial) which is denoted by $P_{\m}(R(\calF); n)$.
The degree of $P_{\m}(R(\calF); n)$ is equal to $\dim R(\calF)$ (see \cite[Chapter~IIIB, Theorem~1]{Ser}).
We can therefore write
$$P_{\m}(R(\calF); n)=e_{\m}(R(\calF))\frac{n^{\dim R(\calF)}}{\dim R(\calF) !}+O(n^{\dim R(\calF)-1}),
$$
where $O(n^t)$ denotes a polynomial of degree less than or equal to $t$ and $e_{\m}(R(\calF))$
is, by definition, the multiplicity of $R(\calF)$ at $\m$ (see \cite[Chapter~VA]{Ser}).  The following result is a special case of \cite[Theorem~14.7]{matsumura}.

\begin{thm}\label{T:mult-sum}
If $\calF$ is a fan of dimension $d$ (i.e.~such that the maximum of the dimension of the cones in $\calF$ is $d$) in $\bbR^b$
then $R(\calF)$  has dimension $d$ and its multiplicity is equal to
$$e_{\m}(R(\calF))=\sum_{\dim \sigma=d} e_{\m}(R(\sigma))
$$
where  $\m$ is the unique graded maximal ideal of the rings in question.
\end{thm}
\begin{proof}
  The theorem is the special case of \cite[Theorem~14.7]{matsumura} where $A=R(\calF)$ and $\mathfrak{q}=\m$.  Indeed, the rings $R(\sigma)$ are the localizations of $R(\calF)$ at minimial primes $\mathfrak{q}$ satisfying $\dim R(\calF)/\mathfrak{q} = d$ by Corollay~\ref{C:min-pr}.
\end{proof}

The above result reduces the computation of the multiplicity of $R(\calF)$ at $\m$ (for a complete fan $\calF$) to that of the affine semigroup rings
$R(\sigma)$ at $\m$, for $\sigma$ a cone of $\calF$ of maximal dimension. These latter multiplicities can be computed geometrically from the affine semigroup $C(\sigma)$, as we now explain following Gelfand--Kapranov--Zelevinsky \cite{GKZ}.

To that aim, we need to recall some definitions. Given a cone $\sigma\in \calF$, set
$C(\sigma)_{\bbZ}:=\langle \sigma \rangle \cap H_{\bbZ}$ and $C(\sigma)_{\bbR}:=\langle \sigma\rangle \cap H_{\bbR}$.
We denote by $\vol_{C(\sigma)}$ the unique translation-invariant measure on  $C(\sigma)_{\bbR}$ such that the
volume of a standard unimodular simplex $\Delta$ (i.e.~$\Delta$ is the convex hull of a basis of $H_{\bbZ}$ together
with  $0$) is $1$. Following \cite[p.184]{GKZ}, denote by $K_+(C(\sigma))$ the convex hull of the set
$C(\sigma)\setminus \{0\}$ and by $K_-(C(\sigma))$ the closure of $\sigma\setminus K_+(C(\sigma))$. The set $K_-(C(\sigma))$ is a bounded (possibly not convex)
lattice polyhedron in $C(\sigma)_{\bbR}$ which is called the {\em subdiagram part} of $C(\sigma)$.

\begin{figure}
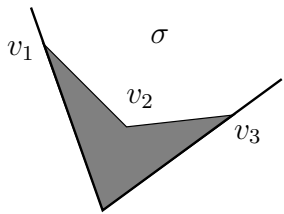

\caption{A two dimensional cone $\sigma$ whose associated semigroup $C(\sigma)$ has Hilbert basis $\calH_{\sigma}=\{v_1,v_2,v_3\}$.
The shaded region is the subdiagram part $K_-(C(\sigma))$ of $C(\sigma)$.}
\label{2Dcone}
\end{figure}

\begin{defi}\cite[Chapter~5, Definition~3.8]{GKZ}
The {\em subdiagram volume}� of $C(\sigma)$ is the natural number
$$u(C(\sigma)):=\vol_{C(\sigma)_{\bbZ}}(K_-(C(\sigma))).$$
\end{defi}

The multiplicity of $R(\sigma)$ at $\m$ can be computed in terms of the subdiagram volume of $C(\sigma)$ as asserted by the following result, whose proof can be found in
\cite[Chapter~5, Theorem~3.14]{GKZ}.

\begin{thm}\label{T:multiGKZ}
The multiplicity of $R(\sigma)$ at $\m$ is equal to
$$e_{\m}(R(\sigma))=u(C(\sigma)).$$
\end{thm}

\section{Geometry of  cographic rings}\label{S:prop-cog}

The aim of this section is to describe the properties of the cographic ring $R(\Gamma)$ associated to a graph $\Gamma$.  The main results are Theorem \ref{T:prop-cog} and the descriptions of the cographic ring in Section~\ref{Ss:desc-cog}.
Recall from Definition \ref{dfnctfr} that $R(\Gamma)$ is the toric face ring associated to the cographic fan $\calF_{\Gamma}^{\perp}$ in $H_1(\Gamma,\bbR)$, which is a polytopal fan by Proposition \ref{Prop: FanCompare}.

According to Proposition \ref{Prop: PosetArr}\eqref{Item: Arr1}, every cone of $\calF_{\Gamma}^{\perp}$ is of the form
$$\sigma(T,\phi):=\bigcap_{e\not\in T}\{(\cdot, \phi(e)) \ge 0\} \bigcap_{e\in T}\{(\cdot, e)=0\},$$
for some uniquely determined element $(T,\phi)\in \OP_{\Gamma}$, i.e.~a totally cyclic orientation $\phi$ on
$\Gamma\setminus T$. We will denote the positive normal affine semigroup associated to $\sigma(T,\phi)$ as in \eqref{E:smgr}
by
$$C(\Gamma\setminus T,\phi):=C(\sigma(T,\phi))=\sigma(Y,\phi)\cap H_1(\Gamma,\bbZ),$$
and its associated affine semigroup ring (as in Definition \ref{D:Aff-smgr}) by
$$
R(\Gamma\setminus T,\phi):=k[C(\Gamma\setminus T,\phi)].
$$

\subsection{Affine semigroup rings $R(\Gamma\setminus T,\phi)$}\label{Ss:aff-cog}

Let us look more closely at the affine semigroup rings $R(\Gamma\setminus T,\phi)$, for a fixed $(T,\phi)\in \OP_{\Gamma}$.

The ring $R(\Gamma\setminus T,\phi)$ is a normal, Cohen--Macaulay domain of dimension equal to $\dim \sigma(T,\phi)=b_1(\Gamma\setminus T)$, as follows from Lemma \ref{L:pro-affsmgr} and Proposition \ref{Prop: PosetArr}\eqref{Item: Arr2}.
However the ring $R(\Gamma\setminus T,\phi)$ need not be Gorenstein, and indeed not even $\bbQ$-Gorenstein, as the following example shows.

\begin{exa}\label{Examp: non-Gor}
Consider the totally cyclic oriented graph $(\Gamma,\phi)$ depicted in Figure \ref{Fig: non-Gor}.

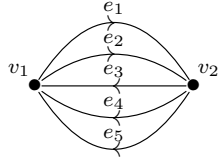
\begin{figure}[ht]
\begin{equation*}
\xymatrix{
 *{\bullet} \ar @{-}@/_1pc/[rr]|-{\SelectTips{cm}{}\object@{<}}^{e_4} \ar @{-} @/_2pc/[rr] |-{\SelectTips{cm}{}\object@{<}}^{e_5} \ar @{-}[rr]|-{\SelectTips{cm}{}\object@{<}}^{e_3}
 \ar@{-} @/^1pc/[rr]|-{\SelectTips{cm}{}\object@{>}}^{e_2}  \ar@{-}@/^2pc/[rr]|-{\SelectTips{cm}{}\object@{>}}^{e_1}^<{v_1}^>{v_2}  &&*{\bullet}
}
\end{equation*}
\hspace{1cm}
\caption{A totally cyclic oriented graph $(\Gamma,\phi)$ with $R(\Gamma,\phi)$ not $\bbQ$-Gorenstein.}\label{Fig: non-Gor}
\end{figure}
Consider the pointed rational polyhedral cone $\sigma(\emptyset,\phi)\subset H_1(\Gamma,\bbR)$
and its dual cone $\sigma(\emptyset,\phi)^{\vee}\subset H_1(\Gamma,\bbR)^{\vee}$  defined by
$$\sigma(\emptyset,\phi)^{\vee}:=\{\ell \in H_1(\Gamma,\bbR)^{\vee}\: :\: \ell(v)\geq 0 \text{ for every } v\in \sigma(\emptyset,\phi) \}.$$
Since for any edge $e\in E(\Gamma)$,
the graph $\Gamma\setminus \{e\}$ with the orientation induced by $\phi$ is totally cyclic, we get that the cone $\sigma(\emptyset,\phi)$ has five
codimension one faces defined  by the equations $\{(\cdot, \phi(e_i))=0\}$ for $i=1,\ldots,5$ (see Corollary \ref{Cor: PosetArr}).
This implies that the extremal rays of $\sigma(\emptyset,\phi)^{\vee}$ are the rays generated by $(\cdot,  \phi(e_i))$ for $i=1,\ldots,5$.

It follows from the proof of \cite[Theorem~3.12]{Dais} that $R(\Gamma,\phi)$ is $\bbQ$-Gorenstein if and only if there exists an element $m\in H_1(\Gamma,\bbQ)$ such that $(m, \phi(e_i))=1$ for every $i=1,\ldots, 5$. However these conditions force $m$ to be equal to $m=\sum_{i=1}^5 \phi(e_i)$, which is a contradiction since
$\partial\left(\sum_{i=1}^5\phi(e_i)\right)=v_1-v_2\neq 0$.

\end{exa}

Denote by $\calH_{(\Gamma\setminus T,\phi)}$ the Hilbert basis (i.e.~the minimal generating set) of the positive affine normal semigroup $C(\Gamma\setminus T,\phi)$. From Lemma \ref{Lemma: DecCyc}, we get the following explicit description of $\calH_{(\Gamma\setminus T,\phi)}$.

\begin{prop}\label{P:Hilb-bas}
The Hilbert basis of $C(\Gamma\setminus T,\phi)$ is equal to
$$\calH_{(\Gamma\setminus T,\phi)}:=\{[\gamma]\: :\: \gamma\in \Cir_{\phi}(\Gamma\setminus T)\}\subset
H_1(\Gamma\setminus T,\bbZ)\subseteq H_1(\Gamma,\bbZ) .$$
\end{prop}

The Hilbert basis $\calH_{(\Gamma\setminus T,\phi)}$ of $C(\Gamma\setminus T,\phi)$ enjoys the following remarkable properties.

\begin{lemma}\label{Lemma: SemiGrp}
Let $(T,\phi)\in \OP_{\Gamma}$. Then
\noindent
	\begin{enumerate}[(i)]
		\item \label{Item: SemiGrp1} The group $\bbZ\cdot \calH_{(\Gamma\setminus T,\phi)}\subseteq
H_1(\Gamma\setminus T ,\bbZ)$  generated by $\calH_{(\Gamma\setminus T,\phi)}$
            coincides with $H_1(\Gamma\setminus T,\bbZ)$.
		\item \label{Item: SemiGrp2} For each $[\gamma]\in \calH_{(\Gamma\setminus T)}$, the ray $\bbR_{\geq 0}\cdot [\gamma]$ is extremal for the cone $\sigma(T,\phi)=\bbR_{\geq 0}\cdot \calH_{(\Gamma\setminus T,\phi)}$.
	\end{enumerate}
\end{lemma}
\begin{proof}
Part (i) follows from Lemma~\ref{Lemma: TotCyclic}(\ref{Item: Gener}).
Part (ii) follows from Proposition \ref{Prop: PosetArr}\eqref{Item: Arr3}.
\end{proof}

We warn the reader that the Hilbert basis $\calH_{(\Gamma\setminus T,\phi)}$ need not  be unimodular
as we show in the Example \ref{Exa: non-unimod} below. Recall that a subset $\calA\subset \bbZ^d$ is said to be {\em unimodular} if $\calA$ spans $\bbR^d$ and, moreover, if we represent the elements of $\calA$ as column vectors of a matrix $A$ with respect to a basis of $\bbZ^d$, then all the nonzero $d\times d$ minors of $A$ have the same absolute value (see \cite[page~70]{Stu}).

 \begin{exa}\label{Exa: non-unimod}
Consider the totally cyclic oriented graph $(\Gamma,\phi)$ depicted in Figure \ref{Fig: non-unimod}.
\begin{figure}[ht]
\begin{equation*}
\xymatrix{
& *{\bullet} \ar@{-}[dl]|-{\SelectTips{cm}{}\object@{<}}^{e^1_0}
\ar @{-}@/_1pc/[dl]|-{\SelectTips{cm}{}\object@{<}}_{e_1^1} \ar@{-}[dr]|-{\SelectTips{cm}{}\object@{>}}_{e_0^2}
\ar @{-}@/^1pc/[dr]|-{\SelectTips{cm}{}\object@{>}}^{e_1^2}& \\
*{\bullet} \ar@{-}[rr]|-{\SelectTips{cm}{}\object@{<}}_{e_0^3}
\ar @{-}@/_1.5pc/[rr]|-{\SelectTips{cm}{}\object@{<}}_{e_1^3}& & *{\bullet}
}
\end{equation*}
\caption{A totally cyclic oriented graph $(\Gamma,\phi)$ with $\calH_{(\Gamma,\phi)}$
not totally unimodular.}
\label{Fig: non-unimod}
\end{figure}
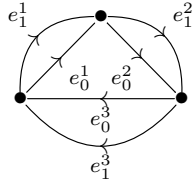
One can check that $b_1(\Gamma)=4$ and that $\calH_{(\Gamma,\phi)}$ consists of the eight elements
$$[\gamma_{ijk}]=\phi(e_i^1)+\phi(e_j^2)+\phi(e_k^3),$$
for $i, j, k \in \{0, 1\}$. The elements
$\calB:=\{[\gamma_{000}], [\gamma_{100}],$ $ [\gamma_{010}],
[\gamma_{001}]\}$ form a basis of $H_1(\Gamma,\bbZ)$. If we order the elements of
$\calH_{(\Gamma,\phi)}$ as
$$\{[\gamma_{000}],  [\gamma_{100}], [\gamma_{010}], [\gamma_{001}],
[\gamma_{110}], [\gamma_{101}], [\gamma_{011}], [\gamma_{111}] \},$$
then the elements of $\calH_{(\Gamma,\phi)}$, with respect to the basis $\calB$, are the column vectors of
the following matrix
$$A=\begin{pmatrix}
1 & 0 & 0 & 0 & -1 & -1& -1& -2 \\
0 & 1 & 0 & 0 & 1 & 1 & 0 & 1 \\
0 & 0 & 1 & 0 & 1 & 0 & 1 & 1 \\
0 & 0 & 0 & 1 & 0 & 1 & 1 & 1\\
\end{pmatrix}.
$$
The minor $A_{1234}$ (i.e.~the minor corresponding to the first four columns) is equal to 1, while
the minor $A_{2348}$ is equal to $2$; hence $\calH_{(\Gamma,\phi)}$ is not unimodular.
\end{exa}

According to \eqref{E:sur-smgrp} and \eqref{E:binom-id}, the affine semigroup ring $R(\Gamma\setminus T,\phi)$ admits the following presentation
\begin{equation}\label{E:pre-cog}
R(\Gamma\setminus T,\phi):=\frac{k[V_{\gamma}\: :\: \gamma\in \Cir_{\phi}(\Gamma\setminus T)]}{I_{(\Gamma\setminus T,\phi)}},
\end{equation}
where $I_{(\Gamma\setminus T,\phi)}:=I_{\sigma(T,\phi)}$ is a binomial ideal, called the toric ideal associated to $\calH_{(\Gamma\setminus T,\phi)}$ in the terminology of \cite[Chapter~4]{Stu}.
The following problem seems interesting.
\begin{prob} \label{Prob: Circuits}
Find generators for the binomial toric ideal $I_{(\Gamma\setminus T,\phi)}$.
\end{prob}

We warn the reader that the toric ideal $I_{(\Gamma\setminus T,\phi)}$ need not to be homogeneous as shown by the following example.

\begin{exa}\label{Exa: non-hom}
Consider the totally cyclicly oriented graph $(\Gamma,\phi)$ depicted in Figure \ref{Fig: non-hom}.
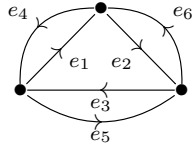
\begin{figure}[ht]
\begin{equation*}
\xymatrix{
& *{\bullet} \ar@{-}[dl]|-{\SelectTips{cm}{}\object@{<}}^{e_1} \ar @{-}@/_1pc/[dl]|-{\SelectTips{cm}{}\object@{>}}_{e_4} \ar@{-}[dr]|-{\SelectTips{cm}{}\object@{>}}_{e_2}
\ar @{-}@/^1pc/[dr]|-{\SelectTips{cm}{}\object@{<}}^{e_6}& \\
*{\bullet} \ar@{-}[rr]|-{\SelectTips{cm}{}\object@{<}}_{e_3} \ar @{-}@/_1pc/[rr]|-{\SelectTips{cm}{}\object@{>}}_{e_5}& & *{\bullet}
}
\end{equation*}
\caption{A totally cyclic oriented graph $(\Gamma,\phi)$ with $I_{(\Gamma,\phi)}$ not homogeneous.}
\label{Fig: non-hom}
\end{figure}

It is easy to see that $b_1(\Gamma)=4$ and that $\calH_{(\Gamma,\phi)}$ consists of the five elements
$$\begin{aligned}
\left[\gamma_1\right] & :=\phi(e_1)+\phi(e_4), \\
[\gamma_2] & :=\phi(e_2)+\phi(e_6), \\
[\gamma_3] & :=\phi(e_3)+\phi(e_5), \\
[\gamma_4] & :=\phi(e_1)+\phi(e_2)+\phi(e_3), \\
[\gamma_5] & :=\phi(e_4)+\phi(e_4)+\phi(e_6). \\
\end{aligned}
$$
The binomial ideal $I_{(\Gamma,\phi)}$ is generated by $V_{\gamma_1}V_{\gamma_2}V_{\gamma_3}-V_{\gamma_4}V_{\gamma_5}$; hence it is not homogeneous.
\end{exa}

\subsection{Descriptions of $R(\Gamma)$ as an inverse limit and as a quotient}\label{Ss:desc-cog}

Using the general results of \S\ref{Ss:2desc}, the ring $R(\Gamma)$ admits two explicit descriptions.

The first description of $R(\Gamma)$ is as an inverse limit of affine semigroup rings (see Proposition \ref{P:pres1}):
\begin{equation}\label{E:invlim}
R(\Gamma)=\varprojlim_{(T,\phi)\in \OP_{\Gamma}} R(\Gamma\setminus T,\phi).
\end{equation}

The second description is a presentation of $R(\Gamma)$ as a quotient of a polynomial ring. In order to make this explicit for $R(\Gamma)$,  observe first that the union of all the Hilbert bases of the cones $\sigma(T,\phi)$, as $(T,\phi)$ varies in $\OP_{\Gamma}$, is equal to the set of all oriented circuits of $\Gamma$, i.e.
\begin{equation}\label{E:tbas-cog}
 \displaystyle \calH_{\calF_{\Gamma}^{\perp}}=\Cyc(\Gamma).
 \end{equation}
Moreover, Corollary \ref{Cor: cyccone} implies that the simplicial complex $\Delta_{\calF^{\perp}}$ introduced in \S\ref{Ss:2desc} coincides with the simplicial complex $\Delta(\Cyc(\Gamma))$ of concordant circuits as in Definition
\ref{Def: Concor}; or in symbols
$$\Delta_{\calF^{\perp}}=\Delta(\Cyc(\Gamma)).$$
From \eqref{E:surR}, Proposition \ref{P:pres2} and Lemma \ref{L:max-orie}, we get the following presentation of $R(\Gamma)$:
\begin{equation}\label{E:pres2-cog}
R(\Gamma)=\frac{ k[V_{\gamma}\: :\: \gamma\in \Cyc(\Gamma)]}{I_{\Gamma}},
\end{equation}
where $I_{\Gamma}:=I_{\calF_{\Gamma}^{\perp}}$ is explicitly given by
\begin{equation}\label{eqncogpres}
I_{\Gamma} =(V_{\gamma}V_{\gamma'}\: :\: \gamma\not\asymp \gamma')+\sum_{(T,\phi)\in \OP_{\Gamma}}
I_{(\Gamma\setminus T,\phi)}
= (V_{\gamma}V_{\gamma'}\: :\: \gamma\not\asymp \gamma')+\sum_{(E(\Gamma)_{\rm sep},\phi)\in \OP_{\Gamma}}
I_{(\Gamma\setminus E(\Gamma)_{\rm sep},\phi)}.
\end{equation}
From Proposition \ref{P:pr-ideals}, we get that the graded prime ideals of $R(\Gamma)$ are given by
\begin{equation}\label{Eqn: Prime-ideals}
{\mathfrak p}_{(T,\phi)}:=(\{X^c\: :\: c\not\in \sigma(T,\phi)\}),
\end{equation}
as $(T,\phi)$ varies in $\OP_{\Gamma}$.

\subsection{Singularities of $R(\Gamma)$}\label{Sec: sing}

In this subsection, we analyze the singularities of the ring $R(\Gamma)$.

\begin{thm}\label{T:prop-cog}
Let $\Gamma$ be a graph and $R(\Gamma)$ its associated cographic ring. Then
\begin{enumerate}[(i)]
\item \label{I:minpr-cog} $R(\Gamma)$ is a reduced finitely generated $k$-algebra of pure dimension equal to $b_1(\Gamma)$.
The minimal prime ideals of $R(\Gamma)$ are
given by $\p_{(E(\Gamma)_{\rm sep},\phi)}$, as $\phi$ varies among all the totally cyclic orientations of $\Gamma\setminus E(\Gamma)_{\rm sep}$.
\item \label{I:Gor-cog} $R(\Gamma)$ is Gorenstein and its canonical module $\omega_{R(\Gamma)}$ is isomorphic to
$R(\Gamma)$ as a graded module.
\item \label{I:semnor-cog} $R(\Gamma)$ is a seminormal ring.
\item \label{I:norm-cog} The normalization of $R(\Gamma)$ is equal to
$$\ov{R(\Gamma)}=\prod_{\phi} R(\Gamma\setminus E(\Gamma)_{\rm sep},\phi),$$
where the product is over all the totally cyclic orientations $\phi$ of $E(\Gamma)\setminus E(\Gamma)_{\rm sep}$.
\item \label{I:slc-cog} The variety $\Spec R(\Gamma)$ is slc.
\item \label{I:emb-cog} The embedded dimension of $R(\Gamma)$ at $\m$ is equal to the cardinality of $\Cyc(\Gamma)$, the set of oriented circuits on $\Gamma$.
\item \label{I:mult-cog} The multiplicity of $R(\Gamma)$ at $\m$ is equal to
$$e_{\m}(R(\Gamma))=\sum_{\phi} e_{\m}(R(\Gamma\setminus E(\Gamma)_{\rm sep},\phi)=
\sum_{\phi} u(C(\Gamma\setminus E(\Gamma)_{\rm sep},\phi))
$$
where the sum is over all the totally cyclic orientations $\phi$ of $\Gamma\setminus E(\Gamma)_{\rm sep}$
and $\m$ is the unique graded maximal ideal of the rings in question.
\end{enumerate}
\end{thm}
\begin{proof}
Part \eqref{I:minpr-cog} follows from Remark \ref{remred}, Corollary \ref{C:min-pr} and Lemma \ref{L:max-orie}.
Part  \eqref{I:Gor-cog} follows Theorem \ref{T:Gore} using that $\calF_{\Gamma}^{\perp}$ is a polytopal fan by Proposition \ref{Prop: FanCompare}.
Part \eqref{I:semnor-cog} follows from Theorem \ref{T:semnorm}\eqref{I:semnorm2}.
Part \eqref{I:norm-cog} follows from Theorem \ref{T:semnorm}\eqref{I:semnorm1} and Lemma \ref{L:max-orie}.
Part \eqref{I:slc-cog} follows from Theorem \ref{T:slc} using that $\calF_{\Gamma}^{\perp}$ is polytopal.
Part \eqref{I:emb-cog} follows from Theorem \ref{T:embdim}\eqref{I:Emb2} and \eqref{E:tbas-cog}.
Part \eqref{I:mult-cog} follows from Theorem \ref{T:mult-sum}, Theorem \ref{T:multiGKZ} and Lemma \ref{L:max-orie}.
\end{proof}

\begin{prob} \label{Prob: Multiplicity}
Express the multiplicity of $R(\Gamma)$ at $\mathfrak m$ in terms of well-known graph invariants.
\end{prob}

\begin{prob} \label{Rmk: sdlt}
Characterize the graphs $\Gamma$ that have the property that $\Spec(R(\Gamma))$ is semi divisorial log terminal.  (See \cite[Definition~1.1]{Fuj1} for the definition of semi divisorial log terminal.)
\end{prob}

Problem~\ref{Rmk: sdlt} is motivated by moduli theory.  The singularities of $R(\Gamma)$ are the singularities that appear on compactified Jacobians, and compactified Jacobian arise as limits of abelian varieties.  In \cite{Fuj2}, Fujino shows that, in a suitable sense, it is possible to degenerate an abelian variety to a semi divisorial log terminal variety. If $R(\Gamma)$ is semi divisorial log terminal, then compactified Jacobians are examples of Fujino's degenerations. For a general discussion of singularities and their role in moduli theory, we direct the reader to \cite{Kol}.

Following the proof of Theorem~\ref{T:slc}, Problem \ref{Rmk: sdlt} is equivalent to the following one: characterize the totally cyclic orientations $\phi$ of a  graph $\Gamma$ that have the property that the pair $(\Spec R(\Gamma,\phi),$ $ -K_{R(\Gamma,\phi)})$  is divisorial log terminal (in the sense of \cite{KM}).
Note that the pair $(\Spec R(\Gamma,\phi),$ $-K_{R(\Gamma,\phi)})$ does not satisfy the stronger condition of being Kawamata log terminal (and so $\Spec R(\Gamma)$ is not semi Kawamata log terminal)  because $-K_{R(\Gamma,\phi)}$ is effective and nonzero.

\section{The cographic ring  $R(\Gamma)$ as a ring of invariants}\label{Sec: presen}

In \cite{CMKV}, the completion of the  ring $R(\Gamma)$ with respect to the maximal ideal $\m=\p_{\un{0}}$  appears  naturally as a ring of invariants.
In this section, we explain this connection.
Consider the multiplicative group
$$\displaystyle T_{\Gamma} := \prod_{v \in V(\Gamma) } \Gm.$$
The elements of $T_{\Gamma}(S)$ for a $k$-scheme $S$ can be written as $\lambda=(\lambda_v)_{v\in V(\Gamma)}$ with
$\lambda_v\in \bbG_m(S)=\calO_S^*$.

Consider the ring
$$A(\Gamma) :=  \frac{k[U_{\el}, U_{\er}\: :\: e\in E(\Gamma)]} {(U_{\el} U_{\er} \colon e \in E(\Gamma))}.$$
If we make the group $T_{\Gamma}$ act on $A(\Gamma)$ via
$$\lambda \cdot U_{\er}= \lambda_{s(\er)}\:U_{\er}\: \lambda_{t(\er)}^{-1},$$
then the invariant subring is described by the following theorem.

\begin{thm}\label{Thm: Pres-inv}
The invariant subring $A(\Gamma)^{T_{\Gamma}}$
is isomorphic to the cographic toric ring $R(\Gamma)$.
\end{thm}

\begin{proof}
We prove the theorem by exhibiting a $k$-basis for the invariant subring that is indexed by  $H_{1}(\Gamma, \bbZ)$ in such a way that multiplication satisfies Eq.~\eqref{Eqn: CographicTimesLaw}.  We argue as follows.  Grade  $A(\Gamma)$ by the $\OC_{1}(\Gamma, \bbZ)$-grading induced by the obvious grading of  $k[U_{\el}, U_{\er}\: :\: e\in E(\Gamma)]$ (so the weight of $U_{\er}$ is $\er$).

This grading is preserved by the  action of $T_{\Gamma}$ on $A(\Gamma)$, so the invariant subring is generated by invariant homogeneous elements.  Furthermore, given a homogeneous  element  $M^{c} = \prod U_{\er}^{a(\er)}$ of weight $c = \sum a(\er) \er$, an element $\lambda \in T_{\Gamma}$ acts as
\begin{align*}
  \lambda \cdot M^{c} =&	\prod_{\er} \lambda_{s(\er)}\:U_{\er}\: \lambda_{t(\er)}^{-1} \\
	    =&	(\prod_{v} \lambda_v^{b(v)}) M^{c},
\end{align*}
where $b(v)$ is defined by $\partial(c) = \sum b(v) v$.  In particular, we see that $M^{c}$ is invariant if and only if $\partial(c)=0$, or in other words $c \in H_{1}(\Gamma, \bbZ)$.

We can conclude that the invariant subring is generated by the homogeneous elements $M^{c}$ whose weight $c$ lies in $H_{1}(\Gamma, \bbZ)$.  In fact, these elements freely generate the invariant subring, because distinct elements have distinct weights.

To complete the proof, observe that  multiplication satisfies
\begin{equation}
	M^c\cdot M^{c'}=
\begin{sis}
	& 0 &\text{ if  $(c,\er)>0$, $(c',\er)<0$ for some $\er$; } \\
	& M^{c+c'} & \text{ otherwise.}
\end{sis}
\end{equation}
The condition that there exists an oriented edge $\er$ with $(c,\er)>0$ and $(c', \er)<0$ is equivalent to the condition that $c$ and $c'$ do not lie in a common  cone, by Corollary~\ref{Cor: cyccone}.  We can conclude that the rule $X^{c} \mapsto M^{c}$ defines an isomorphism between the cographic ring $R(\Gamma)$ and the invariant subring of $A(\Gamma)$.
\end{proof}

\section{A Torelli-type result for $R(\Gamma)$}\label{Sec: Tor}

In this section, we investigate when two graphs give rise to the same cographic toric face ring.
Before stating the result, we need to briefly recall some operations in graph theory introduced in \cite[Section~2]{CV1}.
Two graphs $\Gamma$ and $\Gamma'$ are said to be \emph{cyclic equivalent} (or $2$-isomorphic) if there exists a bijection
$\epsilon: E(\Gamma)\to E(\Gamma')$ inducing a bijection on the circuits. The cyclic equivalence class of $\Gamma$ is denoted
by $[\Gamma]_{\rm cyc}$.
Given a graph $\Gamma$, a \emph{$3$-edge connectivization} of $\Gamma$ is a graph which is obtained from $\Gamma$ by contracting all the separating edges
of $\Gamma$ and by contracting, for every separating pair of edges, one of the two edges.
While a $3$-edge connectivization  of $\Gamma$ is not unique (because of the freedom that we have in performing the second operation), its
cyclic equivalence class is well-defined; it is called the \emph{$3$-edge connected class} of $\Gamma$ and denoted by $[\Gamma]^3_{\rm cyc}$.

\begin{thm}\label{Thm: Tor}
Let $\Gamma$ and $\Gamma'$ be two graphs. Then $R(\Gamma)\cong R(\Gamma')$ if and only if $[\Gamma]_{\rm Cyc}^3=[\Gamma']_{\rm Cyc}^3$.
\end{thm}
\begin{proof}
Assume first that $[\Gamma]_{\rm Cyc}^3=[\Gamma']_{\rm Cyc}^3$. From the proof of \cite[Proposition~3.2.3]{CV1}, it follows that
$\calC_{\Gamma}^{\perp}\cong \calC_{\Gamma'}^{\perp}$, i.e.~that there exists an $\bbR$-linear isomorphism $\phi:H_1(\Gamma,\bbR)\to H_1(\Gamma',\bbR)$
that sends $H_1(\Gamma,\bbZ)$ isomorphically onto $H_1(\Gamma',\bbZ)$ and such that $\phi$ sends the hyperplanes of $\calC_{\Gamma}^{\perp}$
bijectively onto the hyperplanes of $\calC_{\Gamma'}^{\perp}$.  Since $\calF_{\Gamma}^{\perp}$ is the fan induced by the arrangement of hyperplanes
$\calC_{\Gamma}^{\perp}$,
the above map $\phi$ will send the cones of $\calF_{\Gamma}^{\perp}$ bijectively
onto the cones of $\calF_{\Gamma'}^{\perp}$. Therefore the map
$$\begin{aligned}
R(\Gamma) & \to R(\Gamma')\\
X^c& \mapsto X^{\phi(c)},\\
\end{aligned}
$$
is an isomorphism of rings.

Conversely, if $R(\Gamma)\cong R(\Gamma')$ then clearly $\Str(R(\Gamma))\cong \Str(R(\Gamma'))$ (see Definition \ref{D:strata}).
By Corollary \ref{C:strata}, we deduce that $\OP_{\Gamma}\cong \OP_{\Gamma'}$, which implies that  $[\Gamma]_{\rm Cyc}^3=[\Gamma']_{\rm Cyc}^3$
by \cite[Theorem~5.3.2]{CV1}.
\end{proof}

\bibliography{Cographic-toricARXIVFINAL}

\end{document}

%% file: subdiagramscale2.pdf_t
\begin{picture}(0,0)%
\includegraphics{subdiagramscale2.pdf}%
\end{picture}%
\setlength{\unitlength}{3947sp}%
\begingroup\makeatletter\ifx\SetFigFont\undefined%
\gdef\SetFigFont#1#2#3#4#5{%
  \reset@font\fontsize{#1}{#2pt}%
  \fontfamily{#3}\fontseries{#4}\fontshape{#5}%
  \selectfont}%
\fi\endgroup%
\begin{picture}(1762,1319)(5686,-4583)
\put(5701,-3586){\makebox(0,0)[lb]{\smash{{\SetFigFont{12}{14.4}{\familydefault}{\mddefault}{\updefault}{\color[rgb]{0,0,0}$v_1$}%
}}}}
\put(7126,-4111){\makebox(0,0)[lb]{\smash{{\SetFigFont{12}{14.4}{\familydefault}{\mddefault}{\updefault}{\color[rgb]{0,0,0}$v_3$}%
}}}}
\put(6601,-3511){\makebox(0,0)[lb]{\smash{{\SetFigFont{12}{14.4}{\familydefault}{\mddefault}{\updefault}{\color[rgb]{0,0,0}$\sigma$}%
}}}}
\put(6451,-3886){\makebox(0,0)[lb]{\smash{{\SetFigFont{12}{14.4}{\familydefault}{\mddefault}{\updefault}{\color[rgb]{0,0,0}$v_2$}%
}}}}
\end{picture}%